\documentclass[pdflatex,sn-mathphys-num]{sn-jnl}

\usepackage{graphicx}%
\usepackage{multirow}%
\usepackage{amsmath,amssymb,amsfonts}%
\usepackage{amsthm}%
\usepackage{mathrsfs}%
\usepackage[title]{appendix}%
\usepackage{xcolor}%
\usepackage{textcomp}%
\usepackage{manyfoot}%
\usepackage{booktabs}%
\usepackage{algorithm}%
\usepackage{algorithmicx}%
\usepackage{algpseudocode}%
\usepackage{listings}%

\usepackage{subcaption}

\theoremstyle{thmstyleone}%
\newtheorem{theorem}{Theorem}
\theoremstyle{thmstyletwo}%
\newtheorem{remark}{Remark}%

\theoremstyle{thmstylethree}%
\newtheorem{lemma}[theorem]{Lemma}
\newtheorem{assumption}{Assumption}
\begin{document}

\title[$q$-HMC]{Hamiltonian Monte Carlo from $q$-deformed phase-space mechanics}


\author[1]{\fnm{Xiaomei} \sur{Yang}}\email{yangxiaomath@swjtu.edu.cn}

\author*[2]{\fnm{Zhiliang} \sur{Deng}}\email{dengzhl@uestc.edu.cn}


\affil[1]{\orgdiv{School of Mathematics}, \orgname{Southwest Jiaotong University}, \orgaddress{\street{Xi'an Road}, \city{Chengdu}, \postcode{611756}, \state{Sichuan}, \country{China}}}

\affil*[2]{\orgdiv{School of Mathematical Science}, \orgname{University of Electronic Science and Technology of China}, \orgaddress{\street{Xiyuan Ave}, \city{Chengdu}, \postcode{611731}, \state{Sichuan}, \country{China}}}

\abstract{
Hamiltonian Monte Carlo (HMC) constructs Markov transitions by combining
Hamiltonian dynamics with a Metropolis correction. This paper develops a
\(q\)-deformed HMC framework motivated by \(q\)-deformed phase-space mechanics.
Starting from a Lagrangian formulation, we derive the corresponding
\(q\)-Hamiltonian equations and show the formal invariance of the associated
\(q\)-symplectic form. This provides a geometric origin for the proposed
\(q\)-HMC dynamics rather than treating it as an ad hoc finite-difference
modification.
To obtain a computable sampler for ordinary target distributions, we realize
the formal \(q\)-differential operators through Jackson-type differences and
construct a Metropolis-adjusted \(q\)-Hamiltonian proposal. The resulting
transition reduces to classical HMC as \(q\to1\), satisfies detailed balance
under the Metropolis correction, and admits sufficient nondegeneracy conditions
for irreducibility of the marginal position chain.
Numerical experiments on benchmark targets, positive-scale black-box targets,
a mixed positive/signed target, and a PDE-based Bayesian inverse problem show
that the \(q\)-Jackson force acts as a relative perturbation mechanism in the
original variables. It avoids the domain violations and scale mismatch that
can occur with raw additive finite differences, while tracking exact-gradient,
adjoint-force, or carefully scaled finite-difference benchmarks when these are
available.
}




\maketitle

\section{Introduction}\label{sec:intr}

Hamiltonian Monte Carlo (HMC) is a powerful Markov chain Monte Carlo method for sampling from probability distributions of the form
\begin{align}\label{target_distri}
    \pi(x)\propto \exp[-U(x)],
\end{align}
where \(U(x)\) is the potential function. By augmenting the parameter \(x\) with an auxiliary momentum variable \(p\), HMC constructs proposals through Hamiltonian dynamics in the extended phase space. Compared with random-walk-type methods, HMC can generate distant proposals while maintaining high acceptance probabilities, and has become a standard tool in Bayesian computation, statistical machine learning, and scientific computing \cite{Betancourt2017, Bishop2006, Bou-Rabee2020, Duane1987, Durmus2020, Livingstone2019, Mackay2003, Neal1995, Neal2011}.

The performance of HMC is closely tied to the geometric structure of the
underlying Hamiltonian system. In the classical setting, the Hamiltonian flow
preserves the symplectic form and the phase-space volume, which provides the
foundation for stable long-distance proposals. In practical computations,
however, the continuous Hamiltonian flow must be approximated by a numerical
integrator, and the resulting proposal must be corrected by a Metropolis step.
This geometric viewpoint has motivated a variety of extensions of HMC,
including modified or shadow Hamiltonian methods
\cite{Akhmatskaya2008, Radivojevic2020, Heide2021}, alternative choices of
kinetic energy \cite{Livingstone2019_b, Lu2017}, non-Euclidean or
Riemannian-manifold metrics \cite{Betancourt2013, Girolami2011}, and other
structure-preserving proposal mechanisms based on symplectic or
volume-preserving discretizations \cite{Afshar2015, Chao2015, Neal2011}.
In this paper, we pursue a different direction: we study a \(q\)-deformed
analogue of Hamiltonian dynamics and use it to construct a corresponding
\(q\)-HMC sampler.

%
%

The motivation comes from \(q\)-deformed mechanics on \(q\)-commutative spaces. In \(q\)-calculus, ordinary derivatives are replaced by Jackson-type finite differences, and the underlying phase-space algebra is modified by \(q\)-commutation relations. Such structures have been studied in the context of \(q\)-deformed quantum mechanics and non-commutative geometry \cite{Bardek2000, Cerchiai1999, Hinterding1999, Iida1992, Wess1990}. In particular, Lavagno et al. \cite{Lavagno2006} formulated a \(q\)-deformed Hamiltonian system through a \(q\)-Poisson bracket. These constructions suggest a natural question: can one build an HMC-type sampling method by replacing the classical Hamiltonian dynamics with a computable \(q\)-deformed Hamiltonian dynamics?

We answer this question by developing a geometric \(q\)-analogue of HMC. Starting from a Lagrangian formulation, we derive the \(q\)-Hamiltonian equations associated with the \(q\)-commutative phase space. We then prove the invariance of the corresponding \(q\)-symplectic form under the formal \(q\)-Hamiltonian flow. This result may be viewed as a \(q\)-analogue of the classical symplectic invariance and Liouville theorem, formulated within the \(q\)-deformed differential calculus. Since the formal operator-valued dynamics is not directly suitable for numerical simulation, we subsequently introduce a computable realization based on Jackson derivatives.
The resulting deterministic proposal mechanism reduces to classical Hamiltonian dynamics as \(q\to 1\), while for \(q\neq 1\) it replaces ordinary derivatives by \(q\)-Jackson finite differences. Because the computable \(q\)-deformed proposal is not, in general, a canonical Hamiltonian flow in the ordinary Euclidean phase space, its validity as an MCMC transition must be established through a Metropolis correction. We therefore construct a Metropolis-corrected \(q\)-HMC algorithm and prove that the resulting Markov transition satisfies the detailed balance condition with respect to the target distribution.

In addition to the general geometric construction, the \(q\)-Jackson force has a useful computational interpretation. It perturbs variables multiplicatively rather than additively. This feature is not restricted to a single class of targets, but it becomes especially transparent for positive variables. If \(s>0\) and \(y=\log s\), then the symmetric multiplicative perturbations \(s\mapsto e^{\delta/2}s\) and \(s\mapsto e^{-\delta/2}s\), with \(\delta=|\log q|\), correspond to centered perturbations \(y\mapsto y\pm\delta/2\) in the log-coordinate. Thus, for positive black-box targets, the \(q\)-Jackson force can be interpreted as a multiplicative implementation of a log-coordinate finite difference. This observation provides one concrete regime in which the \(q\)-analogue proposal exhibits a visible practical advantage.

The numerical experiments are organized to examine both the validity and the computational behavior of the proposed method. We compare \(q\)-HMC with classical HMC, raw additive finite-difference HMC, and log-coordinate finite-difference HMC. The experiments show that the proposed method is a valid HMC-type sampler and that, in positive-scale and multiscale black-box targets, the \(q\)-Jackson force can avoid the instability caused by raw additive perturbations. In these cases, \(q\)-HMC closely tracks log-coordinate finite-difference HMC and the exact-gradient benchmark, while raw additive finite differences may produce large force errors, large Hamiltonian errors, and poor acceptance rates.


The rest of the paper is organized as follows. Section~\ref{sec:q-analog} derives the
\(q\)-deformed Hamiltonian system and proves the invariance of the associated
\(q\)-symplectic form. Subsequently, a computable \(q\)-Hamiltonian system is presented via the \(q\)-Jackson derivative. 
Section~\ref{sec3:mh} introduces the \(q\)-HMC algorithm
and establishes its validity as a Markov chain Monte Carlo method, including
detailed balance and irreducibility conditions. Section~\ref{sec:numerics} presents the
computational interpretation and numerical behavior of the \(q\)-Jackson force.
In particular, it explains its relation to finite-difference approximations and
demonstrates its scale-consistency advantage for positive, multiscale, and
mixed-coordinate black-box targets. Section~\ref{sec:conclude} concludes the paper.

\section{A \(q\)-Analogue of Hamiltonian Dynamics}\label{sec:q-analog}

%
%

This section develops the geometric \(q\)-analogue of Hamiltonian dynamics that will later serve as the basis for the proposed sampler. We first recall the \(q\)-covariant differential structure on the quantum plane, derive the corresponding \(q\)-Hamiltonian equations from a Lagrangian formulation, and prove the invariance of the associated \(q\)-symplectic form. The same \(q\)-Hamiltonian system was introduced in \cite{Lavagno2006} through a \(q\)-Poisson bracket; the present derivation provides a complementary Lagrangian viewpoint and makes the underlying symplectic structure explicit. These results are formulated in the formal \(q\)-deformed phase space and should be understood as geometric properties of the \(q\)-Hamiltonian system itself, not as a direct invariance statement for a Euclidean target probability measure. We then pass to a computationally tractable realization based on \(q\)-Jackson derivatives, yielding the deterministic \(q\)-deformed vector field used to generate the proposal dynamics. 

\subsection{\(q\)-calculus and \(q\)-Hamiltonian equations}

Let \(\hat{x}\) and \(\hat{p}\) be non-commuting coordinate generators satisfying
the \(q\)-commutation relation
\begin{align}\label{rel1.1}
\hat{p}\hat{x}=q\hat{x}\hat{p},
\qquad
[\hat{p},\hat{x}]_q:=\hat{p}\hat{x}-q\hat{x}\hat{p}=0,
\end{align}
where \(q>0\) is a deformation parameter, with \(q\neq1\) in the deformed case. The classical commutative
case is recovered in the limit \(q\to 1\). This relation is invariant under
the action of the quantum group \(\mathrm{GL}_q(2)\). We denote by
\[
\hat{\mathcal A}=\mathbb C\langle \hat{x},\hat{p}\rangle/
\langle \hat{p}\hat{x}-q\hat{x}\hat{p}\rangle
\]
the coordinate algebra of the two-dimensional quantum plane. Equivalently,
\(\hat{\mathcal A}\) is generated by the non-commuting coordinates
\(\hat{x}\) and \(\hat{p}\), subject to the single relation
\(\hat{p}\hat{x}=q\hat{x}\hat{p}\).
The \(q\)-tangent space \(T\hat{\mathcal A}\) is generated by the
\(q\)-derivatives \(\hat{\partial}_{x}\) and \(\hat{\partial}_{p}\). As
differential operators, their action on the generators is defined by
\begin{align}\label{eq:q-deriv-action}
\hat{\partial}_{x}(\hat{x})=1,\,
\hat{\partial}_{x}(\hat{p})=0,\,
\hat{\partial}_{p}(\hat{x})=0,\,
\hat{\partial}_{p}(\hat{p})=1,
\end{align}
where the parentheses indicate operator action on functions. These operators
satisfy the following \(q\)-commutation relations with the coordinates
\cite{Ubriaco1992}:
\begin{align}\label{eq:q-deriv-algebra}
\begin{aligned}
\hat{\partial}_{p}\hat{x} &= q\hat{x}\hat{\partial}_{p},\\
\hat{\partial}_{x}\hat{p} &= q\hat{p}\hat{\partial}_{x},\\
\hat{\partial}_{p}\hat{p}
&=1+q^2\hat{p}\hat{\partial}_{p}
+(q^2-1)\hat{x}\hat{\partial}_{x},\\
\hat{\partial}_{x}\hat{x}
&=1+q^2\hat{x}\hat{\partial}_{x},\\
\hat{\partial}_{p}\hat{\partial}_{x}
&=q^{-1}\hat{\partial}_{x}\hat{\partial}_{p}.
\end{aligned}
\end{align}
Equation~\eqref{eq:q-deriv-algebra} specifies the algebraic relations between
operators, whereas \eqref{eq:q-deriv-action} specifies their action on the
coordinate generators.

To develop the \(q\)-deformed Hamiltonian formalism, we introduce the module
of one-forms generated by \(d\hat{x}\) and \(d\hat{p}\). Following the
\(\mathrm{GL}_q(2)\)-covariant differential calculus
\cite{Wess1990}, the coordinates and differentials satisfy
\cite{Lavagno2006}
\begin{align}\label{eq:q-diff-forms1}
\begin{aligned}
\hat{x}\,d\hat{p} &= q\,d\hat{p}\,\hat{x},\\
\hat{p}\,d\hat{x} &= q\,d\hat{x}\, \hat{p}+(q^2-1)\,d\hat{p}\,\hat{x},\\
\hat{x}\,d\hat{x} &= q^2\,d\hat{x}\,\hat{x},\\
\hat{p}\,d\hat{p} &= q^2\,d\hat{p}\,\hat{p}.
\end{aligned}
\end{align}
These relations ensure the compatibility of the exterior differential operator
\(d\) with the algebraic structure of the quantum plane.

Throughout this section, we use left \(q\)-derivatives and write one-forms with
the basis differentials on the left and coefficients on the right. Thus, for a
function \(\hat{F}\in\hat{\mathcal A}\),
\begin{align}\label{eq:left-diff-convention1}
d\hat{F}=d\hat{x}\,(\hat{\partial}_{\hat{x}}\hat{F})^L+d\hat{p}\,(\hat{\partial}_{\hat{p}}\hat{F})^L .
\end{align}
Here the superscript (L) denotes left action of the derivative operator on a
function; it does not mean that the resulting coefficient is placed on the left
of the one-form basis. This convention will also be used below for the
variables \((\hat{x},\hat{v})\).

Let \(\hat{v}:=\dot{\hat{x}}\) denote the \(q\)-deformed velocity, where the dot
denotes differentiation with respect to a central time parameter \(t\). Thus
\(t\) commutes with \(\hat{x}\) and \(\hat{p}\). We assume that the time
evolution is compatible with the \(q\)-differential calculus, in the sense that
the exchange relations between coordinates and differentials are preserved
under the time derivative. In particular, the velocity satisfies
\begin{align}\label{eq:v-dp-commutation}
\hat{v}\,d\hat{p}=q\,d\hat{p}\,\hat{v}.
\end{align}

Let the \(q\)-deformed Lagrangian be
\(\hat{\mathcal L}=\hat{\mathcal L}(\hat{x},\hat{v})\). The corresponding
\(q\)-deformed momentum is defined by the left action of the velocity
derivative on the Lagrangian:
\begin{align}\label{qmom}
\hat{p}=q^{1/2}\bigl(\hat{\partial}_{\hat{v}}\hat{\mathcal L}\bigr)^L .
\end{align}
The momentum coordinate is assumed to satisfy the same \(q\)-commutation
relation \eqref{rel1.1} with \(\hat{x}\). The \(q\)-deformed Euler--Lagrange
equation is
\begin{align}\label{eleq}
\frac{d}{dt}\bigl(\hat{\partial}_{\hat{v}}\hat{\mathcal L}\bigr)^L
-\bigl(\hat{\partial}_{\hat{x}}\hat{\mathcal L}\bigr)^L
=0,
\end{align}
which follows from the formal stationary condition for the action
\[
S=\int \hat{\mathcal L}(\hat{x},\hat{v})\,dt .
\]
Combining \eqref{qmom} and \eqref{eleq} gives
\begin{align}\label{pdot}
\dot{\hat{p}}=q^{1/2}\bigl(\hat{\partial}_{\hat{x}}\hat{\mathcal L}\bigr)^L .
\end{align}
Using the convention \eqref{eq:left-diff-convention1}, the differential of the
Lagrangian is
\begin{align}\label{Ld}
\begin{aligned}
d\hat{\mathcal L}
&=d\hat{x}\,\bigl(\hat{\partial}_{\hat{x}}\hat{\mathcal L}\bigr)^L+
d\hat{v}\,\bigl(\hat{\partial}_{\hat{v}}\hat{\mathcal L}\bigr)^L\\
&=q^{-1/2}\left(d\hat{x}\,\dot{\hat{p}}+d\hat{v}\,\hat{p}\right)\\
&=q^{-1/2}\left(d\hat{x}\,\dot{\hat{p}}+d(\hat{v}\,\hat{p})-\hat{v}\,d\hat{p}\right).
\end{aligned}
\end{align}
The \(q\)-deformed Hamiltonian is defined by the Legendre transform
\begin{align}\label{qham_def}
\hat{H}(\hat{x},\hat{p}):=
q^{-1/2}\hat{v}\,\hat{p}-\hat{\mathcal L}(\hat{x},\hat{v}).
\end{align}
Then, using \eqref{Ld}, we have
\begin{align}\label{dH_pre}
d\hat{H}=d\left(q^{-1/2}\hat{v}\hat{p}-\hat{\mathcal L}\right)=q^{-1/2}\left(\hat{v}\,d\hat{p}-d\hat{x}\,\dot{\hat{p}}\right).
\end{align}
By the compatibility relation \eqref{eq:v-dp-commutation},
\(\hat{v}\,d\hat{p}=q\,d\hat{p}\,\hat{v}\), and hence it holds that
\begin{align}\label{dH_ordered}
d\hat{H}=d\hat{x}\,\bigl(-q^{-1/2}\dot{\hat{p}}\bigr)+d\hat{p}\,\bigl(q^{1/2}\hat{v}\bigr).
\end{align}
Comparing \eqref{dH_ordered} with
\[
d\hat{H}=d\hat{x}\,\bigl(\hat{\partial}_{\hat{x}}\hat{H}\bigr)^L+d\hat{p}\,\bigl(\hat{\partial}_{\hat{p}}\hat{H}\bigr)^L ,
\]
we obtain the \(q\)-deformed Hamiltonian equations
\begin{align}\label{Hamilton}
\begin{aligned}
\dot{\hat{x}}
&=
q^{-1/2}\bigl(\hat{\partial}_{\hat{p}}\hat{H}\bigr)^L,\\
\dot{\hat{p}}
&=
-q^{1/2}\bigl(\hat{\partial}_{\hat{x}}\hat{H}\bigr)^L.
\end{aligned}
\end{align}
These equations are consistent with those in \cite{Lavagno2006} and reduce to
the classical Hamiltonian system as \(q\to 1\).

\subsection{Formal invariance of the \(q\)-symplectic form}
\label{subsec:q-symplectic-invariance}

We now formulate the symplectic property of the formal \(q\)-Hamiltonian
system.  Associated with the \(q\)-calculus on the quantum plane, we use the
\(q\)-deformed exterior algebra
\[
\Omega(\hat{\mathcal A})
=
\bigoplus_{k=0}^{2}\Omega^k(\hat{\mathcal A}),
\]
where \(\Omega^k(\hat{\mathcal A})\) denotes the space of \(k\)-forms.  Since
the quantum plane is two-dimensional, \(\Omega^k(\hat{\mathcal A})=0\) for
\(k\geq 3\).  Thus \(\Omega^0(\hat{\mathcal A})=\hat{\mathcal A}\),
\(\Omega^1(\hat{\mathcal A})\) is generated by
\(\{d\hat{x},d\hat{p}\}\), and \(\Omega^2(\hat{\mathcal A})\) is generated by
\(d\hat{x}\wedge d\hat{p}\).  The \(q\)-exterior product satisfies
\cite{Lavagno2006, Wess1990}
\begin{align}
\label{q_wedge}
d\hat{p}\wedge d\hat{x}
=
-q^{-1}d\hat{x}\wedge d\hat{p}.
\end{align}
The \(q\)-symplectic form is defined by
\begin{align}
\label{q_omega_def}
\hat{\omega}
=
q^{-1/2}d\hat{x}\wedge d\hat{p}.
\end{align}
This is the formal \(q\)-deformed analogue of the canonical symplectic form.

Let
\[
\hat{X}_{\hat{H}}
=
\dot{\hat{x}}\hat{\partial}_{\hat{x}}
+
\dot{\hat{p}}\hat{\partial}_{\hat{p}}
\]
be the \(q\)-Hamiltonian vector field determined by the component equations
\eqref{Hamilton}.  With the normalization \eqref{q_omega_def} and the same
left-differential convention used in the Lagrangian derivation, these component
equations are equivalently summarized by the one-form identity
\begin{align}
\label{q_hamiltonian_identity}
i_{\hat{X}_{\hat{H}}}\hat{\omega}
=
d\hat{H},
\end{align}
where \(i_{\hat{X}_{\hat{H}}}\) denotes the \(q\)-deformed interior product.
Indeed, using the \(q\)-contraction rules together with the ordering convention
for one-forms, the left-hand side of \eqref{q_hamiltonian_identity} has the
ordered form
\[
i_{\hat{X}_{\hat{H}}}\hat{\omega}
=
d\hat{x}\,\bigl(-q^{-1/2}\dot{\hat{p}}\bigr)
+
d\hat{p}\,\bigl(q^{1/2}\dot{\hat{x}}\bigr),
\]
which coincides with the expression for \(d\hat H\) obtained in
\eqref{dH_ordered}.  Thus no additional scalar factor multiplying \(d\hat H\)
is introduced; the \(q\)-factors are already contained in the \(q\)-wedge
relation, the normalization of \(\hat\omega\), and the contraction rules.

We next define the infinitesimal action of a formal vector field on
\(q\)-forms.  Let \(\phi_t\) denote the formal \(q\)-Hamiltonian flow generated
by \(\hat{X}_{\hat{H}}\).  On the \(q\)-differential calculus described above,
we define the formal Lie derivative by the Cartan homotopy formula
\begin{align}
\label{Cartan}
\mathcal L_{\hat{X}}\hat{\eta}
:=
d\bigl(i_{\hat{X}}\hat{\eta}\bigr)
+
i_{\hat{X}}\bigl(d\hat{\eta}\bigr),
\qquad
\hat{\eta}\in\Omega(\hat{\mathcal A}).
\end{align}
Here the \(q\)-exterior derivative and the \(q\)-interior product are understood
in the sense of the \(q\)-commutative differential calculus.  Therefore
\eqref{Cartan} is used as the formal definition of the Lie derivative on
\(q\)-forms.  It should not be interpreted as a statement about an ordinary
commutative flow on Euclidean phase space.

\begin{theorem}
\label{thm2.1}
The formal \(q\)-Hamiltonian vector field generated by \(\hat H\) preserves the
\(q\)-symplectic form in the infinitesimal sense
\begin{align}
\mathcal L_{\hat{X}_{\hat{H}}}\hat{\omega}
=
0 .
\end{align}
Equivalently, the formal \(q\)-Hamiltonian flow \(\phi_t\), when interpreted
within the \(q\)-deformed differential calculus, preserves \(\hat{\omega}\).
\end{theorem}

\begin{proof}
By the definition \eqref{Cartan} of the formal Lie derivative,
\begin{align}
\mathcal L_{\hat{X}_{\hat{H}}}\hat{\omega}
=
d\bigl(i_{\hat{X}_{\hat{H}}}\hat{\omega}\bigr)
+
i_{\hat{X}_{\hat{H}}}\bigl(d\hat{\omega}\bigr).
\end{align}
Since \(\hat{\omega}\in\Omega^2(\hat{\mathcal A})\) is a top-degree form on the
two-dimensional quantum plane, \(d\hat{\omega}=0\).  Hence
\begin{align}
\mathcal L_{\hat{X}_{\hat{H}}}\hat{\omega}
=
d\bigl(i_{\hat{X}_{\hat{H}}}\hat{\omega}\bigr).
\end{align}
Using the \(q\)-Hamiltonian one-form identity
\eqref{q_hamiltonian_identity}, we obtain
\begin{align}
\mathcal L_{\hat{X}_{\hat{H}}}\hat{\omega}
&=
d\bigl(d\hat{H}\bigr)  \nonumber\\
&=
d^2\hat{H}
=
0 .
\end{align}
The last equality follows from the nilpotency of the \(q\)-exterior derivative.
This proves the formal invariance of the \(q\)-symplectic form.
\end{proof}

\begin{remark}
Theorem~\ref{thm2.1} is a formal \(q\)-deformed analogue of the classical
symplectic invariance of Hamiltonian flow.  When \(q=1\), the \(q\)-calculus
reduces to the ordinary exterior calculus, \eqref{Cartan} becomes the standard
Cartan formula, and the theorem reduces to the usual preservation of the
canonical symplectic form.

For \(q\neq1\), however, the statement is made inside the non-commutative
\(q\)-differential calculus.  It should not be identified directly with
Lebesgue-volume preservation of the computable Euclidean realization introduced
below.  The validity of the Markov transition kernel generated by that
realization will be established separately through the Metropolis--Hastings
correction.
\end{remark}

\subsection{A computable realization through Jackson derivatives}

Although the formal \(q\)-Hamiltonian flow preserves the \(q\)-symplectic form,
the operator-valued variables \(\hat{x}\) and \(\hat{p}\) are not directly
suitable for numerical simulation. We therefore use a standard realization of
the \(q\)-algebra and its differential calculus through dilation and Jackson
derivative operators \cite{Lavagno2006, Ubriaco1992}. This realization provides
a deterministic vector field on the classical Euclidean phase space. Since the
resulting map is not a canonical Hamiltonian flow in the ordinary sense, its
use in a Metropolis sampler requires a separate correction, which will be
discussed in Section~\ref{sec3:mh}.

The realization is given by the replacements
\begin{align}\label{q_realization}
\hat{x}\to x,\,\,
\hat{p}\to pD_x,\,\,
\hat{\partial}_{\hat{x}}\to \mathcal D_x,\,\,
\hat{\partial}_{\hat{p}}\to \mathcal D_pD_x,
\end{align}
where \(D_x\) is the dilation operator in the \(x\)-direction,
\[
D_x f(x, p)=f(qx, p),
\]
and \(\mathcal D_x\), \(\mathcal D_p\) are the forward Jackson derivatives
\begin{align}\label{jackson_derivatives}
\begin{aligned}
\mathcal D_x f(x,p)
&=
\frac{f(q^2x,p)-f(x,p)}{(q^2-1)x},\\
\mathcal D_p f(x,p)
&=
\frac{f(x,q^2p)-f(x,p)}{(q^2-1)p}.
\end{aligned}
\end{align}
These are the standard Jackson derivatives used in the realization of the
\(q\)-differential calculus. They are one-sided multiplicative finite
differences and reduce to the ordinary partial derivatives as \(q\to1\), provided
the corresponding derivatives exist. As usual in \(q\)-calculus, the values at \(x=0\) or \(p=0\) are understood
through the corresponding limiting values whenever these limits exist. The
numerical implementation near such points will be specified in Section~\ref{sec3:mh}.
\begin{remark}
The Jackson derivatives in \eqref{jackson_derivatives} are the standard
one-sided Jackson derivatives associated with the formal \(q\)-differential
calculus. For numerical force approximation, however, we also use a symmetric
log-normalized version.  For a scalar variable \(x\neq0\), set
\(\delta=|\log q|\) and define
\begin{align}\label{sym_jackson_derivative}
\mathcal D^{\mathrm{sym}}_{\delta,x} f(x)
=
\frac{f(e^{\delta/2}x)-f(e^{-\delta/2}x)}{\delta x}.
\end{align}
Equivalently,
\[
    x\,\mathcal D^{\mathrm{sym}}_{\delta,x}f(x)
    =
    \frac{f(e^{\delta/2}x)-f(e^{-\delta/2}x)}{\delta},
\]
which is a centered finite difference with respect to the log-coordinate
\(\log |x|\) when \(x\neq0\).  This expression recovers the ordinary derivative
as \(\delta\to0\), whenever the ordinary derivative exists.  The standard
one-sided form is used in the formal \(q\)-Hamiltonian derivation, while the
symmetric log-normalized form is used in the finite-difference force
approximations in the numerical experiments. In the numerical sections, unless
otherwise stated, the term \(q\)-Jackson force refers to this symmetric
log-normalized form.
\end{remark}

Applying the realization \eqref{q_realization} to the formal equations
\eqref{Hamilton} gives the computable \(q\)-Hamiltonian system
\begin{align}\label{q_ham}
\begin{aligned}
\dot{x}
&=q^{-1/2}\mathcal D_p[H(qx, p)],\\
\dot{p}
&=-q^{1/2}\mathcal D_x[H(x, p)],
\end{aligned}
\end{align}
where \(x\), \(p\), and \(H\) now denote classical variables and a classical
Hamiltonian function. This system reduces to the classical Hamiltonian
equations in the limit \(q\to 1\), provided the ordinary derivatives exist.

Following the terminology of \(q\)-deformed mechanics, we define the bracket-like operation
\begin{align}\label{poisson_bracket}
\{f, g\}_q
=q^{-1/2}\mathcal D_p D_x g\cdot \mathcal D_x f
-q^{1/2}\mathcal D_x g\cdot \mathcal D_p f.
\end{align}
Then \eqref{q_ham} can be written in the compact form
\begin{align}\label{q_anal}
\dot{x}=\{x, H\}_q,\qquad
\dot{p}=\{p, H\}_q.
\end{align}
We call \eqref{q_ham}, or equivalently \eqref{q_anal}, the computable
\(q\)-analogue of Hamiltonian dynamics. This formulation provides the
deterministic \(q\)-deformed vector field used to generate the sampling
proposal in the next section.

\section{A Metropolis-adjusted \(q\)-Hamiltonian sampler}\label{sec3:mh}

In this section, we use the computable \(q\)-Hamiltonian system
\eqref{q_ham} to construct a Metropolis-adjusted sampling method. After the
Jackson realization introduced in Section~\ref{sec:q-analog}, the variables
\(x\) and \(p\) are treated as classical Euclidean variables. Unless otherwise
specified, volume preservation, Jacobian matrices, and determinants in this
section are understood with respect to the Lebesgue measure on
\(\mathbb R^{2d}\).

The section is organized as follows. We first define the \(q\)-Hamiltonian
proposal map and the corresponding Metropolis-adjusted algorithm. We then prove
that the resulting transition satisfies detailed balance with respect to the
extended target density. Finally, we record a sufficient nondegeneracy condition for
irreducibility of the marginal position chain.

\subsection{The \(q\)-Hamiltonian proposal}

Let the target density \(\pi(x)\) be defined on an open set
\(\mathcal X\subseteq\mathbb R^d\) with respect to the Lebesgue measure and satisfy
\begin{align}
    \pi(x)\propto \exp[-U(x)],\qquad x\in\mathcal X.
\end{align}
We set \(U(x)=+\infty\) outside \(\mathcal X\). In the full-support case,
\(\mathcal X=\mathbb R^d\).
We introduce an auxiliary momentum variable \(p\in\mathbb R^d\) and consider
the extended density
\begin{align}\label{extended_target}
    \rho(x,p)
    \propto
    \exp[-H(x,p)],\,\,
    H(x,p)=U(x)+K(p),
\end{align}
where the kinetic energy \(K\) is assumed to be even:
\begin{align}\label{even_kinetic}
    K(-p)=K(p).
\end{align}
The separable form in \eqref{extended_target} is used throughout this section.
For a general nonseparable Hamiltonian \(H(x,p)=U(x)+K(x,p)\), the explicit
leapfrog splitting below is not generally reversible without using implicit
generalized leapfrog substeps.
Let \(G_q(x)\) and \(v_q(p)\) denote the potential-gradient approximation and velocity induced by the
computable \(q\)-Hamiltonian system:
\begin{align}\label{q_force_velocity}
    G_q(x)=q^{1/2}\mathcal D_x U(x),
    \,\,
    v_q(p)=q^{-1/2}\mathcal D_p K(p).
\end{align}
Here \(\mathcal D_x\) and \(\mathcal D_p\) are applied componentwise. At points
where a denominator in a Jackson derivative vanishes, the corresponding
limiting value is used whenever it exists. Since \(K\) is even, the velocity
\(v_q\) is odd:
\[
    v_q(-p)=-v_q(p).
\]
As \(q\to1\), \(G_q(x)\to\nabla U(x)\) and
\(v_q(p)\to\nabla K(p)\), provided that the ordinary derivatives exist.
In the theoretical construction, \(G_q\) denotes the force induced by the chosen
Jackson realization. The detailed-balance proof below only requires that
\(G_q(x)\) depends on \(x\) alone and that the leapfrog substeps are explicit
shears. In the numerical experiments, we use the symmetric log-normalized
Jackson force defined in \eqref{sym_jackson_derivative}.

For the separable Hamiltonian considered here, the computable \(q\)-Hamiltonian
system takes the form
\begin{align}
    \dot{x}=v_q(p),
    \qquad
    \dot{p}=-G_q(x).
\end{align}
For a step size \(h>0\), define the two explicit shear maps
\begin{align}\label{q_shear_maps}
\begin{aligned}
    A_h(x,p)
    &=
    \bigl(x,p-hG_q(x)\bigr),\\
    B_h(x,p)
    &=
    \bigl(x+h v_q(p),p\bigr).
\end{aligned}
\end{align}
Combining these two updates in symmetric order gives one \(q\)-leapfrog step:
\begin{align}\label{q_leapfrog_map}
    \Phi_{q,h}=A_{h/2}\circ B_h\circ A_{h/2}.
\end{align}
For \(N_{\mathrm{lf}}\) leapfrog steps, let
\begin{align}\label{q_multistep_map}
    \Phi_{q,h,N_{\mathrm{lf}}}
    =
    (\Phi_{q,h})^{N_{\mathrm{lf}}}.
\end{align}
Let \(T(x,p)=(x,-p)\) be the momentum-flip map and define
\begin{align}\label{proposal_map}
    \Psi_{q,h,N_{\mathrm{lf}}}
    =
    T\circ \Phi_{q,h,N_{\mathrm{lf}}}.
\end{align}
The momentum flip is included for the reversibility proof. Since the position
component is unchanged by \(T\), it has no effect on the accepted value of \(x\).

Starting from the current position \(x_{\mathrm{curr}}\), we first refresh the
momentum from the kinetic density \(\nu(p)\propto\exp[-K(p)]\), then apply the
\(q\)-leapfrog map followed by the momentum flip. The resulting proposal is
accepted or rejected by a Metropolis rule. The complete transition is
summarized in Algorithm~\ref{alg:qhmc_revised}.

\begin{algorithm}
\caption{Metropolis-adjusted \(q\)-Hamiltonian sampler}
\label{alg:qhmc_revised}
\begin{algorithmic}[1]
\Require Target density \(\pi(x)\propto \exp[-U(x)]\) on \(\mathcal X\), current state \(x_{\mathrm{curr}}\), step size \(h\), number of leapfrog steps \(N_{\mathrm{lf}}\), deformation parameter \(q\), kinetic energy \(K\).
\Ensure Next state \(x_{\mathrm{next}}\).

\State Sample \(p\sim \nu\), where \(\nu(p)\propto \exp[-K(p)]\).
\State Set \(z=(x_{\mathrm{curr}},p)\).
\State Compute \(\tilde z=\Phi_{q,h,N_{\mathrm{lf}}}(z)\) by the \(q\)-leapfrog map \eqref{q_multistep_map}.
\State Set \(z^*=\Psi_{q,h,N_{\mathrm{lf}}}(z)=T\tilde z\).
\State Compute the acceptance probability
\[
\alpha(z, z^*)
=\min\left\{1,\exp[-H(z^*)+H(z)]\right\}.
\]
\State Draw \(u\sim \mathrm{Uniform}(0,1)\).
\If{\(u<\alpha(z,z^*)\)}
    \State Set \(x_{\mathrm{next}}=x^*\), the position component of \(z^*\).
\Else
    \State Set \(x_{\mathrm{next}}=x_{\mathrm{curr}}\).
\EndIf
\State \Return \(x_{\mathrm{next}}\).
\end{algorithmic}
\end{algorithm}

The rest of this section establishes the validity of Algorithm~\ref{alg:qhmc_revised}. The proof follows the standard Metropolis-adjusted deterministic-proposal argument: we first show that the $q$-leapfrog proposal is reversible under momentum flip and preserves the Lebesgue volume on the extended phase space. These two properties imply that the usual Metropolis acceptance probability in Algorithm~\ref{alg:qhmc_revised} yields detailed balance with respect to 
$\rho(x,p)\propto \exp[-H(x,p)]$.

\subsection{Detailed balance and invariance} The purpose of this subsection is to establish that the Metropolis-adjusted \(q\)-Hamiltonian proposal leaves the extended target density \(\rho(x,p)\propto \exp[-H(x,p)]\) invariant. The proof follows the standard Metropolis-adjusted deterministic-proposal argument. We first verify that the proposal map is reversible under momentum flip and preserves the Euclidean volume. These two properties imply that the usual Metropolis acceptance probability yields detailed balance. 

\begin{lemma}
\label{lem:q_leapfrog}
For the separable Hamiltonian used in Algorithm~\ref{alg:qhmc_revised}, assume that $K(-p)=K(p)$. Then the map $\Phi_{q,h,N_{\mathrm{lf}}}$ satisfies
\begin{align}
T\circ \Phi_{q,h,N_{\mathrm{lf}}}\circ T
=
\Phi_{q,h,N_{\mathrm{lf}}}^{-1}.
\end{align}
Moreover, $\Phi_{q,h,N_{\mathrm{lf}}}$, and hence $\Psi_{q,h,N_{\mathrm{lf}}}=T\circ\Phi_{q,h,N_{\mathrm{lf}}}$, is volume-preserving with respect to the Lebesgue measure on $\mathbb R^{2d}$.
\end{lemma}

\begin{proof}
Because $G_q(x)$ is independent of $p$, the momentum update $A_h$ is an explicit shear and satisfies
\begin{align}
A_h^{-1}=A_{-h}.
\end{align}
Because $v_q(p)$ is independent of $x$, the position update $B_h$ is also an explicit shear and satisfies
\begin{align}
B_h^{-1}=B_{-h}.
\end{align}
Furthermore,
\begin{align}
T\circ A_h\circ T=A_{-h},
\qquad
T\circ B_h\circ T=B_{-h},
\end{align}
where the second identity uses $v_q(-p)=-v_q(p)$. Therefore,
\begin{align}
T\circ \Phi_{q,h}\circ T
&=
(T\circ A_{h/2}\circ T)
\circ
(T\circ B_h\circ T)
\circ
(T\circ A_{h/2}\circ T)
\nonumber\\
&=
A_{-h/2}\circ B_{-h}\circ A_{-h/2}
=
\Phi_{q,h}^{-1}.
\end{align}
Iterating this identity gives
\begin{align}
T\circ \Phi_{q,h,N_{\mathrm{lf}}}\circ T
=
\Phi_{q,h,N_{\mathrm{lf}}}^{-1}.
\end{align}

By the definitions in \eqref{q_shear_maps}, $A_h$ changes only the momentum component while keeping the position fixed, and $B_h$ changes only the position component while keeping the momentum fixed. Hence both maps are shear maps on the Euclidean phase space. More explicitly, the Jacobian matrix of $A_h$ is block lower triangular with identity matrices on the diagonal, whereas the Jacobian matrix of $B_h$ is block upper triangular with identity matrices on the diagonal. Therefore,
\begin{align}
|\det D A_h(x,p)|
= |\det D B_h(x,p)|=1.
\end{align}
It follows that the composition $\Phi_{q,h,N_{\mathrm{lf}}}$ preserves the Lebesgue volume. Since the momentum flip $T$ also has absolute Jacobian determinant one, the same holds for $\Psi_{q,h,N_{\mathrm{lf}}}$. This completes the proof.
\end{proof}

\begin{theorem}
\label{thm:qhmc_db}
The transition kernel on the extended phase space defined by Algorithm~\ref{alg:qhmc_revised} satisfies detailed balance with respect to
\begin{align}
\rho(x,p)\propto \exp[-H(x,p)].
\end{align}
Consequently, the marginal position chain leaves $\pi(x)\propto \exp[-U(x)]$ invariant.
\end{theorem}

\begin{proof}
By Lemma~\ref{lem:q_leapfrog}, the proposal map $\Psi_{q,h,N_{\mathrm{lf}}}$ is an involution:
\begin{align}
\Psi_{q,h,N_{\mathrm{lf}}}^{-1}
=
\Psi_{q,h,N_{\mathrm{lf}}},
\end{align}
and it preserves the Lebesgue volume. Therefore, the deterministic proposal has the same change-of-variables structure as in standard Metropolis-adjusted reversible deterministic proposals.
Let $z^*=\Psi_{q,h,N_{\mathrm{lf}}}(z)$. Since
\[
\alpha(z,z^*)
=
\min\left\{1,\frac{\rho(z^*)}{\rho(z)}\right\},
\]
we have
\begin{align}
\rho(z)\alpha(z,z^*)
&=
\min\{\rho(z),\rho(z^*)\}
\nonumber\\
&=
\rho(z^*)\alpha(z^*,z).
\end{align}
This proves detailed balance for the accept--reject transition on the extended phase space. The momentum refreshment step samples $p$ from $\nu(p)\propto \exp[-K(p)]$ and therefore preserves the extended density $\rho(x,p)$. Hence the composition of momentum refreshment and the Metropolis-adjusted $q$-Hamiltonian proposal leaves $\rho$ invariant. Marginalizing $\rho(x,p)$ over $p$ gives the target density $\pi(x)$.
\end{proof}

\begin{remark}
No $q$-Jacobian determinant is used in the Metropolis acceptance probability above. This is deliberate. After the Jackson realization, the Markov chain is defined on the Euclidean phase space, and the target density $\rho(x,p)$ is defined with respect to the Lebesgue measure. Hence the relevant Jacobian in the Metropolis--Hastings correction is the Jacobian of the actual proposal map. For the separable $q$-leapfrog map \eqref{q_leapfrog_map}, this determinant is one because the component updates are triangular shears. The $q$-symplectic structure discussed in Section~\ref{sec:q-analog} motivates the proposal, while the validity of the Markov kernel is established here through the standard Metropolis argument.
\end{remark}

\subsection{Irreducibility}
\label{subsec:irreducibility}

Detailed balance proves invariance of the target distribution, but it does not
by itself guarantee that the chain can explore the support of the target.  We
therefore record a simple sufficient condition for irreducibility of the
marginal position chain.  The result is not intended to provide a convergence
rate.  Its role is only to state an accessibility condition for the
Metropolis-adjusted \(q\)-Hamiltonian sampler.

Let
\begin{align}\label{Gx_def}
G_x(p)
=
\operatorname{pr}_x\left(\Psi_{q,h,N_{\mathrm{lf}}}(x,p)\right)
\end{align}
denote the position component of the proposal generated from the initial state
\((x,p)\), where \(\operatorname{pr}_x\) is the projection onto the position
variable.  If \(P\sim \nu\), where \(\nu(p)\propto \exp[-K(p)]\) is the momentum
density, then \(G_x(P)\) is the proposal position generated from the current
state \(x\).

\begin{assumption}
\label{assumption:nondegenerate}
For every \(x\in\mathcal X\), the pushforward measure
\[
    Q_x(A)
    :=
    \int_{\{p:\,G_x(p)\in A\}}\nu(p)\,dp,
    \qquad A\in\mathcal B(\mathcal X),
\]
is equivalent to the Lebesgue measure restricted to \(\mathcal X\). Equivalently,
the proposal position \(G_x(P)\) admits a density \(q_x(x')\) satisfying
\[
    q_x(x')>0
    \quad
    \text{for Lebesgue-a.e. }x'\in\mathcal X .
\]
\end{assumption}

\begin{theorem}
\label{thm:irreducibility}
Assume that \(U\) is continuous and finite on \(\mathcal X\), so that
\(\pi(x)\propto\exp[-U(x)]\) is strictly positive on \(\mathcal X\), and that
\(U(x)=+\infty\) outside \(\mathcal X\). Assume also that the momentum density
\(\nu(p)\propto\exp[-K(p)]\) is strictly positive on \(\mathbb R^d\). If
Assumption~\ref{assumption:nondegenerate} holds, then the position chain
generated by Algorithm~\ref{alg:qhmc_revised} is Lebesgue-irreducible on
\(\mathcal X\).
\end{theorem}

\begin{proof}
Fix \(x\in\mathcal X\), and let \(A\subset\mathcal X\) be a Borel set with
positive Lebesgue measure.  By Assumption~\ref{assumption:nondegenerate},
\[
    Q_x(A)
    =
    \int_{\{p:\,G_x(p)\in A\}}\nu(p)\,dp
    >0 .
\]
Thus the set of momenta that generate proposed positions in \(A\) has positive
\(\nu\)-probability.

For any such momentum \(p\), write
\[
    z=(x,p),
    \qquad
    z'=\Psi_{q,h,N_{\mathrm{lf}}}(x,p).
\]
Since \(G_x(p)\in A\subset\mathcal X\), the proposed position lies in the
support of the target.  Because \(U\) and \(K\) are finite on the states under
consideration, the Metropolis acceptance probability
\[
    \alpha(z,z')
    =
    \min\{1,\exp[-H(z')+H(z)]\}
\]
is strictly positive.  Hence
\[
    P(x,A)
    \ge
    \int_{\{p:\,G_x(p)\in A\}}
    \alpha\!\left((x,p),\Psi_{q,h,N_{\mathrm{lf}}}(x,p)\right)
    \nu(p)\,dp
    >0.
\]
Therefore the position chain can reach every Borel subset of \(\mathcal X\)
with positive Lebesgue measure with positive probability. This proves
Lebesgue-irreducibility on \(\mathcal X\).
\end{proof}

\begin{remark}
Assumption~\ref{assumption:nondegenerate} is a condition on the projected
proposal distribution, not on the full \(q\)-Hamiltonian map. This distinction
is important. Even if the full phase-space map
\(\Psi_{q,h,N_{\mathrm{lf}}}\) is reversible or locally invertible, the
projected map \(p\mapsto G_x(p)\) need not be a global diffeomorphism. The
assumption only requires that the momentum refreshment induces a nondegenerate
proposal distribution on the position space.
\end{remark}

\begin{remark}
For one \(q\)-leapfrog step and a full-support target, the nondegeneracy
condition can be checked directly in common cases. In this case,
\[
    G_x(p)
    =
    x+h\,v_q\left(p-\frac{h}{2}G_q(x)\right).
\]
Thus the projected proposal is nondegenerate whenever the \(q\)-velocity map
\(v_q\) is a diffeomorphism whose image has full Lebesgue support. For example,
for a diagonal quadratic kinetic energy
\[
    K(p)=\frac12 p^T M^{-1}p,
\]
with \(M\) positive definite and diagonal, the corresponding \(q\)-velocity is
an invertible linear map. Hence the proposal position has a strictly positive
density on \(\mathbb R^d\), and Assumption~\ref{assumption:nondegenerate} is
satisfied in the full-support case.
\end{remark}

\begin{remark}
A full geometric-ergodicity result would require additional drift and
minorization conditions. Such results are delicate even for standard HMC,
where the tail behavior of the target distribution plays a decisive role. In
the present work, we restrict ourselves to invariance and irreducibility and
leave quantitative convergence bounds for future work.
\end{remark}



\section{Numerical tests}\label{sec:numerics}

We now examine the numerical behavior of the proposed \(q\)-HMC method.  The
purpose of this section is not to claim that \(q\)-HMC outperforms exact-gradient
HMC when an analytical gradient is available.  Instead, we focus on black-box
settings in which the potential can be evaluated pointwise but the gradient is
unavailable.  Particular attention is given to positive-scale targets, where
the multiplicative structure of the \(q\)-Jackson perturbation becomes
especially visible.

\subsection{General setting, competing methods, and diagnostics}
\label{subsec:numerical-setup}

A recurring setting in the following experiments is a positive variable
\[
    s>0,\qquad y=\log s .
\]
The Markov chain is run in the log-coordinate \(y\), while the potential is
queried as a black-box function \(U_s(s)\) of the original positive variable.
For a scalar positive variable, the corresponding log-coordinate potential is
\[
    V(y)=U_s(e^y)-y ,
\]
where the term \(-y\) is the Jacobian correction.  We write
\[
    G(y)=V'(y)=sU_s'(s)-1,\qquad s=e^y,
\]
for the potential gradient in the log-coordinate.  The HMC momentum equation is
therefore \(\dot p=-G(y)\).  For a vector
\(s=(s_1,\ldots,s_d)\in\mathbb R_+^d\), the same convention is used
componentwise:
\[
    V(y)=U_s(e^y)-\sum_{i=1}^d y_i,\qquad s_i=e^{y_i}.
\]
In black-box problems, \(U_s'(s)\) or \(\nabla U_s(s)\) is not available and the
potential gradient must be approximated from function evaluations.

We compare the following force constructions.

\emph{Exact-gradient HMC.}  This method uses the analytical potential gradient
\(G(y)=\nabla V(y)\).  It is included only as an oracle benchmark and is not
treated as a black-box method.

\emph{Raw additive finite-difference HMC.}  This method approximates the
derivative by perturbing the original positive variable additively.  In the
scalar case,
\[
    U_s'(s)\approx
    \frac{U_s(s+h_s)-U_s(s-h_s)}{2h_s},
\]
with a one-sided modification when \(s-h_s\le 0\).  The resulting
log-coordinate potential-gradient approximation is
\[
    \widehat G_{\mathrm{add}}(y)
    =
    s\frac{U_s(s+h_s)-U_s(s-h_s)}{2h_s}-1.
\]
This method uses a fixed absolute perturbation \(h_s\), so its effective
relative perturbation is \(h_s/s\).  It can therefore become nonlocal when
\(s\) is small.

\emph{Log-coordinate finite-difference HMC.}  This method applies a centered
finite difference directly in the log-coordinate:
\[
    \widehat G_{\log}(y)
    =
    \frac{V(y+h_y)-V(y-h_y)}{2h_y}.
\]
It serves as a scale-aware finite-difference reference.

\emph{Symmetric \(q\)-Jackson HMC.}  This method perturbs the original positive
variable multiplicatively.  With \(\delta=|\log q|\), it uses
\[
    \widehat G_q(y)
    =
    \frac{U_s(e^{\delta/2}s)-U_s(e^{-\delta/2}s)}{\delta}-1.
\]
Since
\[
    e^{\delta/2}s=e^{y+\delta/2},\qquad
    e^{-\delta/2}s=e^{y-\delta/2},
\]
the symmetric \(q\)-Jackson force can be viewed, in this positive-scale setting,
as a centered finite difference in the log-coordinate implemented through
black-box evaluations of \(U_s(s)\).  Its perturbation size is controlled by the
relative deformation level \(\delta=|\log q|\), rather than by an absolute
additive step.

The following diagnostics are reported.  ``Acc.'' denotes the average
Metropolis acceptance rate after burn-in.  ``KS'' denotes the coordinate-wise
Kolmogorov--Smirnov distance between the empirical marginal distribution and a
reference marginal distribution; for multidimensional examples, the mean over
coordinates is reported.  ``ESS'' denotes the effective sample size, averaged
over coordinates when \(d>1\).  When the acceptance rate is extremely small, the
chain can be nearly stuck and standard ESS estimates may become misleading; in
such cases ESS is not reported.  ``Force err.'' denotes the root mean squared
potential-gradient approximation error,
\[
    \left[
    \mathbb E_{\mathrm{ref}}
    \left\{
    \frac{1}{d}\|\widehat G(Y)-G(Y)\|^2
    \right\}
    \right]^{1/2},
\]
where the expectation is approximated using samples from a reference target.
Finally, we define the median positive Hamiltonian error by
\[
    \mathrm{med.}(\Delta H_+)
    :=
    \operatorname{median}\{\Delta H_j:\Delta H_j>0\},
    \qquad
    \Delta H_j=H(z_j')-H(z_j),
\]
computed over proposed HMC trajectories.  Large positive Hamiltonian errors
indicate poor energy conservation and usually lead to low Metropolis acceptance.

\subsection{One-dimensional positive-scale black-box targets}
\label{ex:positive-scale-blackbox}
In this first example, we test four one-dimensional positive-scale targets.  This experiment
serves as a mechanism test: it isolates the effect of additive versus
multiplicative force approximations when the potential is queried only through
the original positive variable \(U_s(s)\).

Cases A and B are lognormal targets.  Specifically, \(Y=\log S\sim
N(\mu,\sigma^2)\), so that the potential with respect to the original variable
\(s\) is
\[
    U_s(s)=\frac{(\log s-\mu)^2}{2\sigma^2}+\log s .
\]
Case A uses \(\mu=-3\), so that the target mass is concentrated near
\(s=\exp(-3)\), representing a small-scale positive variable.  Case B uses
\(\mu=2\), so that the target is concentrated near \(s=\exp(2)\), representing
a large-scale positive variable.  
Case C is a two-scale lognormal mixture.  Equivalently, the log-scale variable
\(Y=\log S\) follows a two-component Gaussian mixture,
\[
    p_Y(y)
    =
    \frac12\phi(y;\mu_1,\sigma_1^2)
    +
    \frac12\phi(y;\mu_2,\sigma_2^2),
\]
where
\[
    \phi(y;\mu,\sigma^2)
    =
    \frac{1}{\sqrt{2\pi\sigma^2}}
    \exp\left\{
        -\frac{(y-\mu)^2}{2\sigma^2}
    \right\}.
\]
The corresponding density of \(S\) is
\[
    p_S(s)=\frac{1}{s}p_Y(\log s),\qquad s>0.
\]
In the experiment, we set \(\mu_1=-3\) and \(\mu_2=2\), so that the two mixture
components represent a small-scale and a large-scale region, respectively. We use \(\sigma_1=\sigma_2=1\) in the numerical experiment.
This case tests whether the force approximation remains stable when the same
positive variable can move between substantially different scales.
%
%
%
%
%
%
%
Case D is a simple Bayesian scale posterior.  Let \(s>0\) be an unknown
standard deviation and suppose
\[
    z_i\mid s\sim N(0,s^2),\qquad i=1,\ldots,n .
\]
The likelihood depends on \(s\) through the sufficient statistic
\[
    S_z=\sum_{i=1}^n z_i^2
\]
and is proportional to
\[
    L(s;z)\propto s^{-n}\exp\left(-\frac{S_z}{2s^2}\right).
\]
We assign a Gaussian prior to the log-scale variable,
\[
    y=\log s,\qquad y\sim N(\mu_0,\sigma_0^2),
\]
or equivalently a lognormal prior to \(s\).  The posterior potential in the
original variable is, up to an additive constant,
\[
    U_s(s)
    =
    \frac{(\log s-\mu_0)^2}{2\sigma_0^2}
    +\log s
    +n\log s
    +\frac{S_z}{2s^2}.
\]
The corresponding log-coordinate potential is
\[
    V(y)=U_s(e^y)-y
    =
    \frac{(y-\mu_0)^2}{2\sigma_0^2}
    +ny
    +\frac{S_z}{2e^{2y}},
\]
and the exact log-coordinate potential gradient is
\[
    G(y)=V'(y)
    =
    \frac{y-\mu_0}{\sigma_0^2}
    +n
    -S_ze^{-2y}.
\]
In the experiment, we set \(n=10\), \(\mu_0=-1\), \(\sigma_0=2\), and
\(S_z=n\exp(2y_\ast)\) with \(y_\ast=-3\).  Thus the likelihood favors a small
scale \(s\approx\exp(-3)\).

Table~\ref{tab:positive-scale-blackbox} summarizes the results.  Cases A and D
show the clearest deterioration of raw additive finite differences because the
target mass is concentrated at a small positive scale.  Case B serves as a
control: when the typical scale is large, the same additive perturbation is
relatively small and raw additive FD behaves similarly to the scale-aware
methods.  Case C shows that the same effect persists in a two-scale mixture.

\begin{table*}[htbp]
\centering
\caption{
Performance comparison on one-dimensional positive-scale black-box targets.
Smaller KS, force error, and \(\mathrm{med.}(\Delta H_+)\), and larger ESS are
preferred. A: small-scale lognormal; B: large-scale lognormal; C:  two-scale lognormal mixture; D: Bayesian scale posterior.
}
\label{tab:positive-scale-blackbox}
\scriptsize
\begin{tabular}{llccccc}
\toprule
Case & Method & Acc. & KS & ESS & Force err. & \(\mathrm{med.}(\Delta H_+)\) \\
\midrule
A
& Exact-gradient HMC              & 0.999 & 0.008 & 2000.0 & \(0\) & \(6.967{\times}10^{-4}\) \\
& Raw additive FD        & 0.842 & 0.015 &  737.6 & \(1.166{\times}10^{0}\) & \(2.101{\times}10^{-1}\) \\
& Log-coordinate FD      & 0.999 & 0.008 & 2000.0 & \(5.957{\times}10^{-15}\) & \(6.967{\times}10^{-4}\) \\
& \(q\)-Jackson, \(q=1.05\) & 0.999 & 0.008 & 2000.0 & \(4.446{\times}10^{-15}\) & \(6.967{\times}10^{-4}\) \\
& \(q\)-Jackson, \(q=0.95\) & 0.999 & 0.008 & 2000.0 & \(1.970{\times}10^{-15}\) & \(6.967{\times}10^{-4}\) \\
\midrule
B
& Exact-gradient HMC              & 0.999 & 0.008 & 2000.0 & \(0\) & \(6.967{\times}10^{-4}\) \\
& Raw additive FD        & 1.000 & 0.008 & 2000.0 & \(3.078{\times}10^{-5}\) & \(6.965{\times}10^{-4}\) \\
& Log-coordinate FD      & 0.999 & 0.008 & 2000.0 & \(4.803{\times}10^{-15}\) & \(6.967{\times}10^{-4}\) \\
& \(q\)-Jackson, \(q=1.05\) & 0.999 & 0.008 & 2000.0 & \(4.060{\times}10^{-15}\) & \(6.967{\times}10^{-4}\) \\
& \(q\)-Jackson, \(q=0.95\) & 0.999 & 0.008 & 2000.0 & \(1.962{\times}10^{-15}\) & \(6.967{\times}10^{-4}\) \\
\midrule
C
& Exact-gradient HMC              & 1.000 & 0.009 & 1541.8 & \(0\) & \(4.875{\times}10^{-4}\) \\
& Raw additive FD        & 0.935 & 0.007 &  974.9 & \(4.141{\times}10^{-1}\) & \(6.676{\times}10^{-2}\) \\
& Log-coordinate FD      & 1.000 & 0.009 & 1541.8 & \(3.978{\times}10^{-7}\) & \(4.875{\times}10^{-4}\) \\
& \(q\)-Jackson, \(q=1.05\) & 1.000 & 0.009 & 1541.8 & \(2.364{\times}10^{-6}\) & \(4.875{\times}10^{-4}\) \\
& \(q\)-Jackson, \(q=0.95\) & 1.000 & 0.009 & 1541.8 & \(2.612{\times}10^{-6}\) & \(4.875{\times}10^{-4}\) \\
\midrule
D
& Exact-gradient HMC              & 0.999 & 0.006 & 2000.0 & \(0\) & \(1.254{\times}10^{-3}\) \\
& Raw additive FD        & 0.615 & 0.062 &  113.9 & \(4.704{\times}10^{0}\) & \(8.440{\times}10^{-1}\) \\
& Log-coordinate FD      & 0.998 & 0.007 & 2000.0 & \(2.539{\times}10^{-3}\) & \(1.478{\times}10^{-3}\) \\
& \(q\)-Jackson, \(q=1.05\) & 0.996 & 0.006 & 2000.0 & \(1.511{\times}10^{-2}\) & \(5.933{\times}10^{-3}\) \\
& \(q\)-Jackson, \(q=0.95\) & 0.995 & 0.006 & 2000.0 & \(1.671{\times}10^{-2}\) & \(6.517{\times}10^{-3}\) \\
\bottomrule
\end{tabular}
\end{table*}

The table shows a consistent scale effect.  In the small-scale cases, raw
additive FD produces much larger force errors and Hamiltonian errors than
log-coordinate FD and \(q\)-Jackson.  In the large-scale lognormal case, the
same absolute additive step is relatively small, and the performance gap
essentially disappears.  Thus the deterioration is not caused by the target
distribution itself, but by the mismatch between additive finite differences
and positive-scale geometry.  The \(q\)-Jackson force behaves similarly to a
log-coordinate finite difference, while remaining expressible through
multiplicative evaluations of the original black-box potential \(U_s(s)\).

\subsection{High-dimensional synthetic positive-scale tests}
\label{subsec:positive-scale-highdim}

We next examine whether the scale effect observed in the one-dimensional
examples persists in higher-dimensional positive-scale targets.  This
experiment is intended as a synthetic stress test.  Its purpose is not to
introduce a new model class, but to check whether the error caused by raw
additive finite differences accumulates as the number of positive variables
increases.

We first consider independent lognormal components,
\[
    y_i=\log s_i\sim N(\mu,\sigma^2),
    \qquad i=1,\ldots,d ,
\]
so that the black-box potential in the original positive variable is
\[
    U_s(s)
    =
    \sum_{i=1}^d
    \left\{
        \frac{(\log s_i-\mu)^2}{2\sigma^2}
        +\log s_i
    \right\}.
\]
The corresponding log-coordinate potential is
\[
    V(y)=U_s(e^y)-\sum_{i=1}^d y_i
    =
    \sum_{i=1}^d
    \frac{(y_i-\mu)^2}{2\sigma^2}.
\]
This deliberately simple target isolates the numerical effect of additive
versus multiplicative force approximation without introducing posterior
correlations or multimodality.

As a second synthetic high-dimensional example, we consider independent
Bayesian scale posteriors.  For each coordinate, let \(s_i>0\) be an unknown
scale parameter and suppose that the corresponding likelihood favors a
small-scale posterior.  This setting is used to test whether additive
finite-difference errors become more severe when several positive scale
parameters are sampled jointly.

We use the same competing force approximations and diagnostics as in
Section~\ref{subsec:numerical-setup}.  The dimension \(d\) is increased while
the positive scale is kept in a small-scale regime.  Table~\ref{tab:positive-scale-highdim}
reports representative results for \(d=2\) and \(d=50\).

\begin{table*}[htbp]
\centering
\caption{
Representative results for high-dimensional positive-scale targets.
``Prod.'' denotes the product target and ``Scale post.'' denotes the Bayesian
scale posterior. Smaller KS distance and \(\mathrm{med.}(\Delta H_+)\), and larger
acceptance rate and ESS are preferred. ESS is not reported when the chain is
nearly stuck due to an extremely low acceptance rate.
}
\label{tab:positive-scale-highdim}
\scriptsize
\setlength{\tabcolsep}{3.5pt}
\renewcommand{\arraystretch}{1.05}
\begin{tabular}{@{}llccccc@{}}
\toprule
Target & Method & Dim. & Acc. & KS & ESS & \(\mathrm{med.}(\Delta H_+)\) \\
\midrule
Prod.
& Exact-gradient HMC        & \(2\)  & 0.999 & 0.015 & 1200.0 & \(1.293{\times}10^{-3}\) \\
& Additive FD      & \(2\)  & 0.614 & 0.025 & 544.9  & \(8.900{\times}10^{-1}\) \\
& Log-FD           & \(2\)  & 0.999 & 0.016 & 1200.0 & \(1.383{\times}10^{-3}\) \\
& \(q\)-Jac., \(q=1.05\) & \(2\)  & 0.999 & 0.014 & 1200.0 & \(1.349{\times}10^{-3}\) \\
& \(q\)-Jac., \(q=0.95\) & \(2\)  & 0.999 & 0.012 & 1200.0 & \(1.337{\times}10^{-3}\) \\
\cmidrule{2-7}
& Exact-gradient HMC        & \(50\) & 0.996 & 0.014 & 1200.0 & \(8.896{\times}10^{-3}\) \\
& Additive FD      & \(50\) & 0.023 & 0.294 & --     & \(1.024{\times}10^{1}\) \\
& Log-FD           & \(50\) & 0.994 & 0.013 & 1200.0 & \(9.152{\times}10^{-3}\) \\
& \(q\)-Jac., \(q=1.05\) & \(50\) & 0.994 & 0.014 & 1200.0 & \(9.001{\times}10^{-3}\) \\
& \(q\)-Jac., \(q=0.95\) & \(50\) & 0.996 & 0.014 & 1200.0 & \(8.969{\times}10^{-3}\) \\
\midrule
Scale post.
& Exact-gradient HMC        & \(2\)  & 0.996 & 0.017 & 1200.0 & \(2.772{\times}10^{-3}\) \\
& Additive FD      & \(2\)  & 0.236 & 0.045 & 205.3  & \(6.540{\times}10^{0}\) \\
& Log-FD           & \(2\)  & 0.997 & 0.015 & 1200.0 & \(3.491{\times}10^{-3}\) \\
& \(q\)-Jac., \(q=1.05\) & \(2\)  & 0.993 & 0.016 & 1200.0 & \(8.604{\times}10^{-3}\) \\
& \(q\)-Jac., \(q=0.95\) & \(2\)  & 0.995 & 0.018 & 1200.0 & \(8.819{\times}10^{-3}\) \\
\cmidrule{2-7}
& Exact-gradient HMC        & \(50\) & 0.986 & 0.015 & 1200.0 & \(2.815{\times}10^{-2}\) \\
& Additive FD      & \(50\) & 0.000 & 0.431 & --     & \(1.031{\times}10^{2}\) \\
& Log-FD           & \(50\) & 0.980 & 0.016 & 1200.0 & \(2.894{\times}10^{-2}\) \\
& \(q\)-Jac., \(q=1.05\) & \(50\) & 0.964 & 0.016 & 1200.0 & \(4.860{\times}10^{-2}\) \\
& \(q\)-Jac., \(q=0.95\) & \(50\) & 0.968 & 0.016 & 1200.0 & \(5.380{\times}10^{-2}\) \\
\bottomrule
\end{tabular}
\end{table*}

The table shows that the scale mismatch of raw additive finite differences
becomes more severe in high dimensions.  For the product target, the raw
additive FD method already has a lower acceptance rate and a much larger
Hamiltonian error at \(d=2\).  When the dimension increases to \(d=50\), its
acceptance rate drops to \(0.023\), the KS distance increases to \(0.294\), and
\(\mathrm{med.}(\Delta H_+)\) reaches \(1.024\times10^{1}\).  By contrast,
exact HMC, log-coordinate FD, and the two \(q\)-Jackson methods remain close to
one another.

The same pattern is stronger for the Bayesian scale posterior.  At \(d=50\),
raw additive FD is essentially unusable: the acceptance rate is numerically
zero, the KS distance increases to \(0.431\), and
\(\mathrm{med.}(\Delta H_+)\) reaches \(1.031\times10^{2}\).  The
log-coordinate and \(q\)-Jackson methods retain high acceptance rates and small
KS distances.  These results show that the scale effect observed in the
one-dimensional examples is not a one-dimensional artifact.  Coordinate-wise
additive finite-difference errors can accumulate across positive variables,
whereas multiplicative perturbations remain stable in relative scale.

These synthetic high-dimensional tests serve as an intermediate step between
the one-dimensional mechanism examples and the Bayesian variance-component
model studied next.

\subsection{Variance-component positive-scale black-box targets}
\label{subsec:variance-component-blackbox}

We next consider a class of positive-scale black-box targets arising from
variance-component inference.  This example is used as the main Bayesian
positive-scale test.  The purpose is again not to outperform an exact analytical
gradient, but to examine whether the multiplicative \(q\)-Jackson perturbation
remains stable when the posterior scale varies over several orders of magnitude
and when the number of positive scale parameters increases.

Let \(\tau_i>0\), \(i=1,\ldots,d\), be variance components and set
\(y_i=\log \tau_i\).  Conditional on \(\tau_i\), observations are generated from
\[
    z_{ij}\mid \tau_i \sim N(0,\tau_i^2),\qquad j=1,\ldots,n_i,
\]
with the log-scale prior
\[
    y_i\sim N(\mu_0,\sigma_0^2).
\]
Writing \(S_i=\sum_{j=1}^{n_i}z_{ij}^2\), the negative log posterior in the
log-scale variable \(y=(y_1,\ldots,y_d)\), up to an additive constant, is
\begin{equation}
\label{eq:variance-component-potential}
    V(y)
    =
    \sum_{i=1}^d
    \left\{
        \frac{(y_i-\mu_0)^2}{2\sigma_0^2}
        + n_i y_i
        + \frac{S_i}{2\exp(2y_i)}
    \right\}.
\end{equation}
The sampler is implemented in \(y\), but the black-box evaluations are made in
the original positive variable \(\tau=\exp(y)\).  The exact-gradient HMC method
uses the analytical log-scale gradient
\[
    \frac{\partial V}{\partial y_i}
    =
    \frac{y_i-\mu_0}{\sigma_0^2}
    + n_i
    - S_i\exp(-2y_i),
\]
and serves only as an oracle benchmark.

We use the same finite-difference force approximations and diagnostics as in
Section~\ref{subsec:numerical-setup}.  In the main \(q\)-Jackson comparison, we
set \(q=\exp(h_y)\) with \(h_y=0.02\), so that \(|\log q|\) is exactly matched
to the log-coordinate finite-difference step.  We also include \(q=1.05\) as a
coarser multiplicative perturbation.

Two diagnostic sweeps are reported.  In the scale sweep, the number of variance
components is fixed and the true log-scale value \(\eta_{\rm true}\) is varied
from very small to moderate scale.  In the dimension sweep,
\(\eta_{\rm true}\) is fixed in a small-scale regime and the number of variance
components is increased.

\begin{figure}[p]
    \centering
    \begin{subfigure}{0.47\textwidth}
        \centering
        \includegraphics[width=\textwidth,height=0.23\textheight,keepaspectratio]{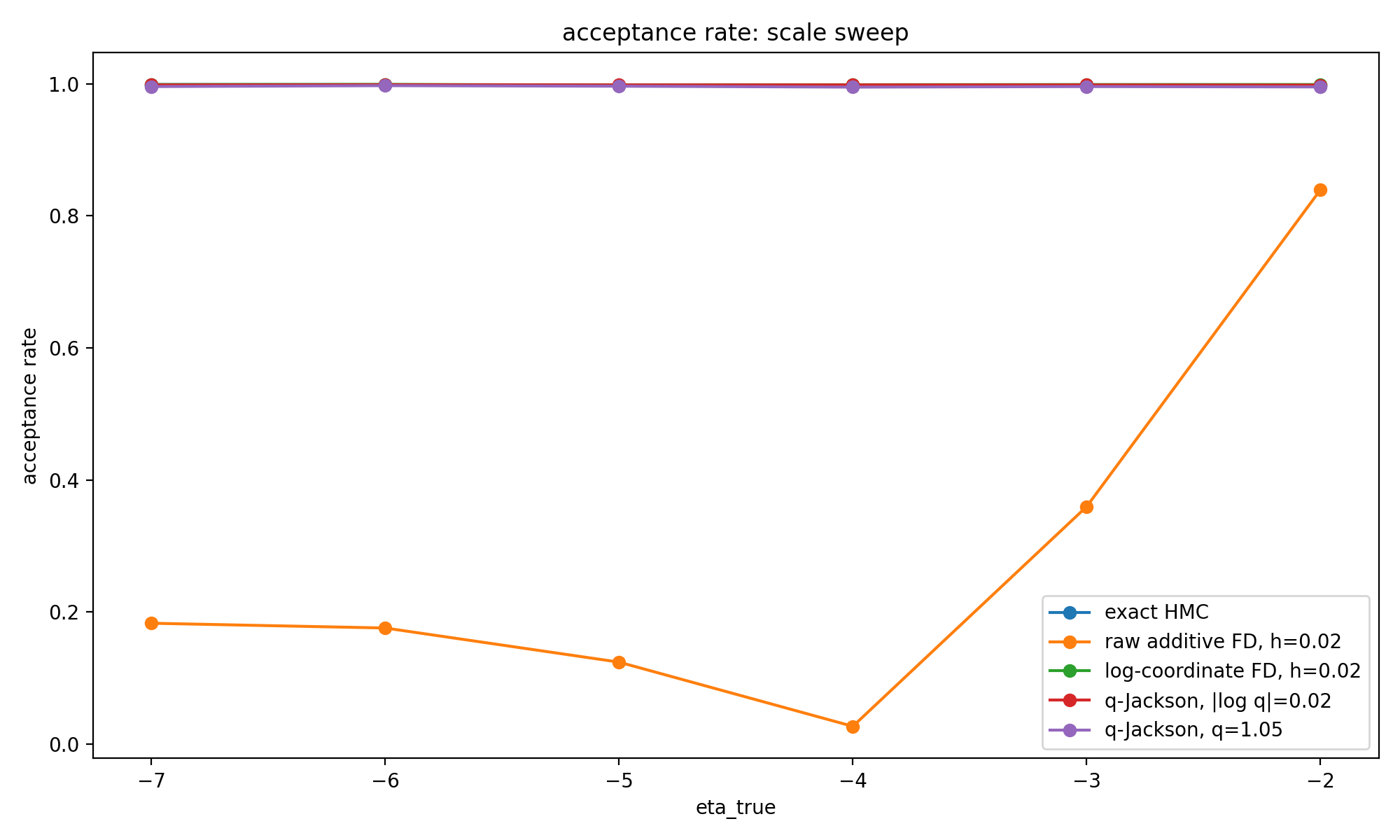}
        \caption{Acceptance rate}
    \end{subfigure}
    \hfill
    \begin{subfigure}{0.47\textwidth}
        \centering
        \includegraphics[width=\textwidth,height=0.23\textheight,keepaspectratio]{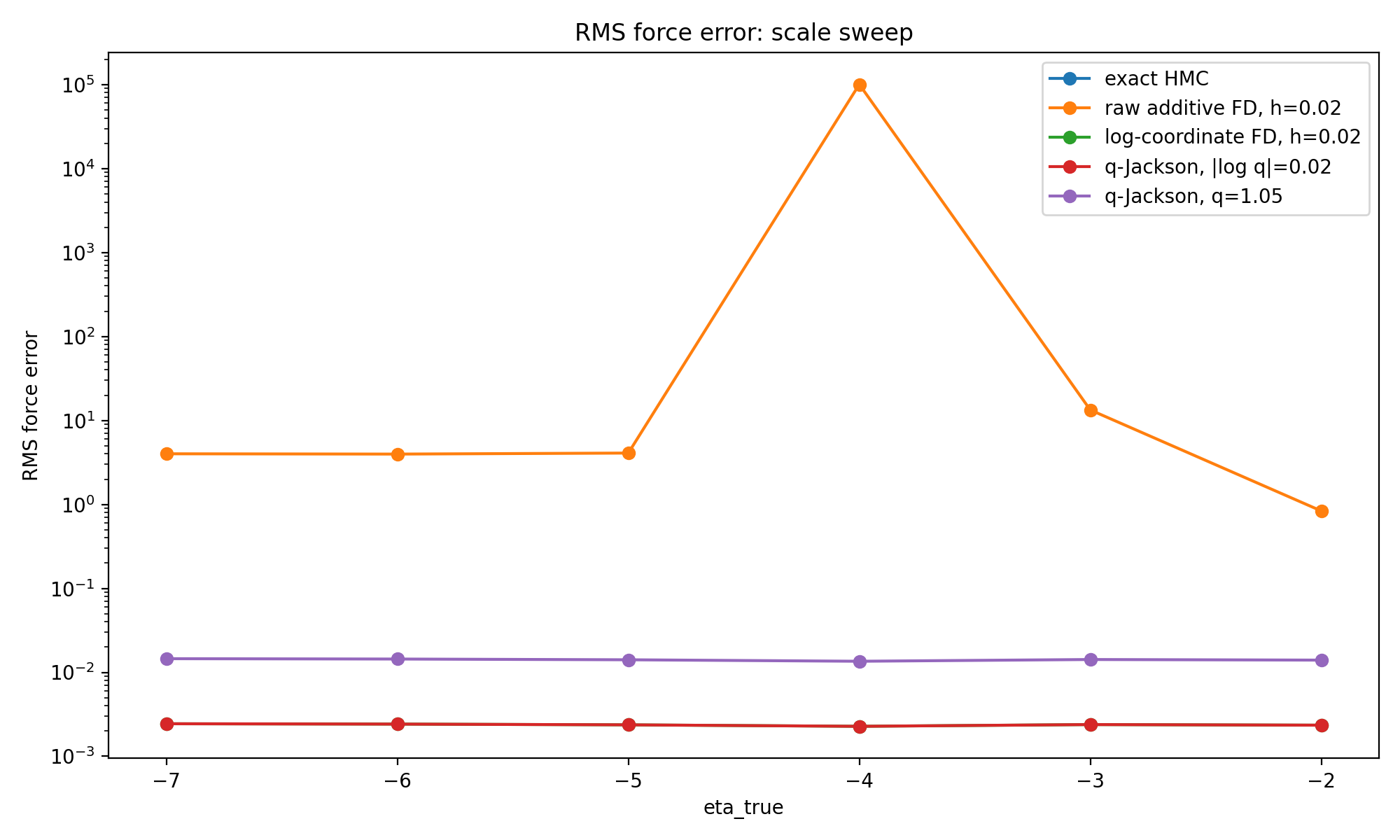}
        \caption{RMS force error}
    \end{subfigure}

    \vspace{0.25em}

    \begin{subfigure}{0.47\textwidth}
        \centering
        \includegraphics[width=\textwidth,height=0.23\textheight,keepaspectratio]{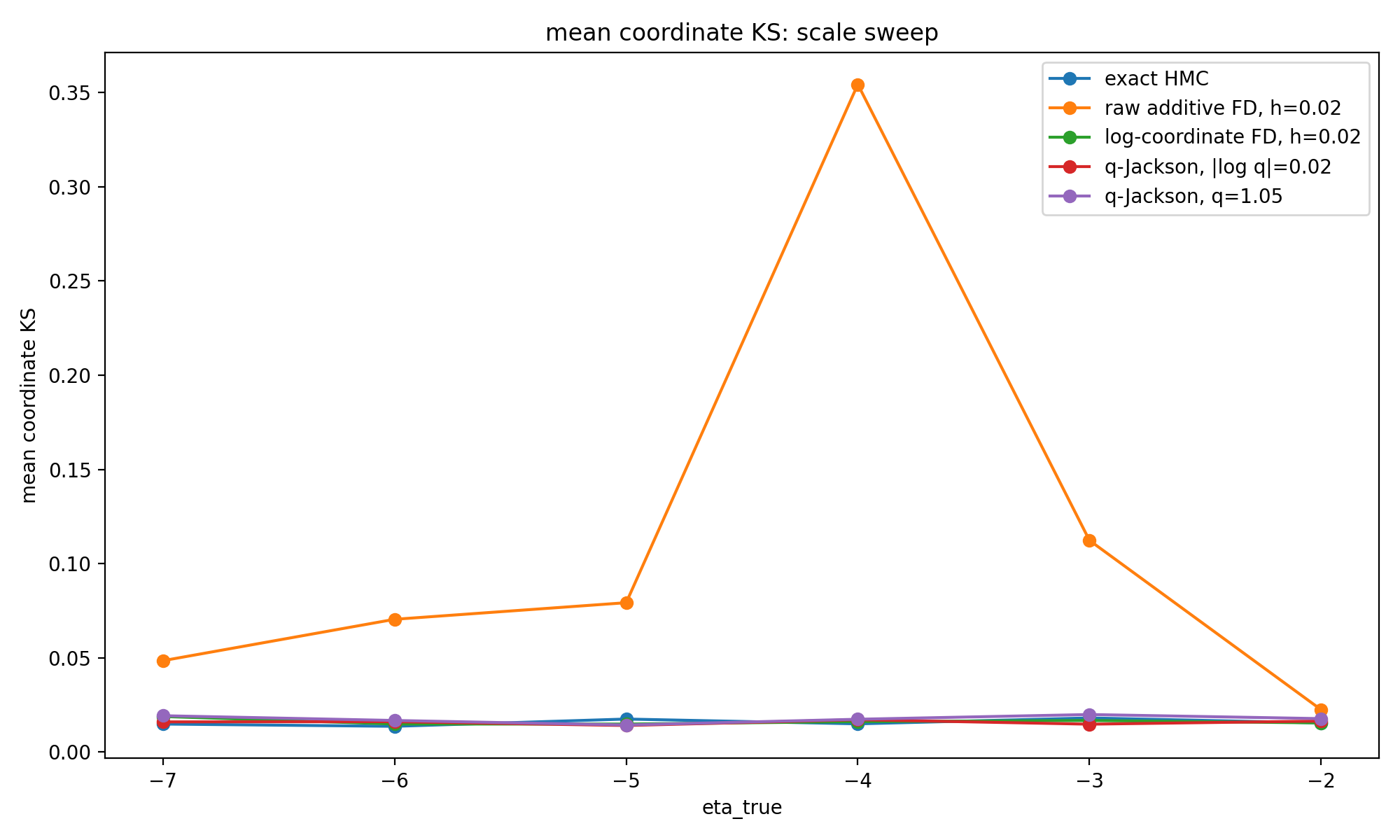}
        \caption{Mean coordinate KS distance}
    \end{subfigure}
    \hfill
    \begin{subfigure}{0.47\textwidth}
        \centering
        \includegraphics[width=\textwidth,height=0.23\textheight,keepaspectratio]{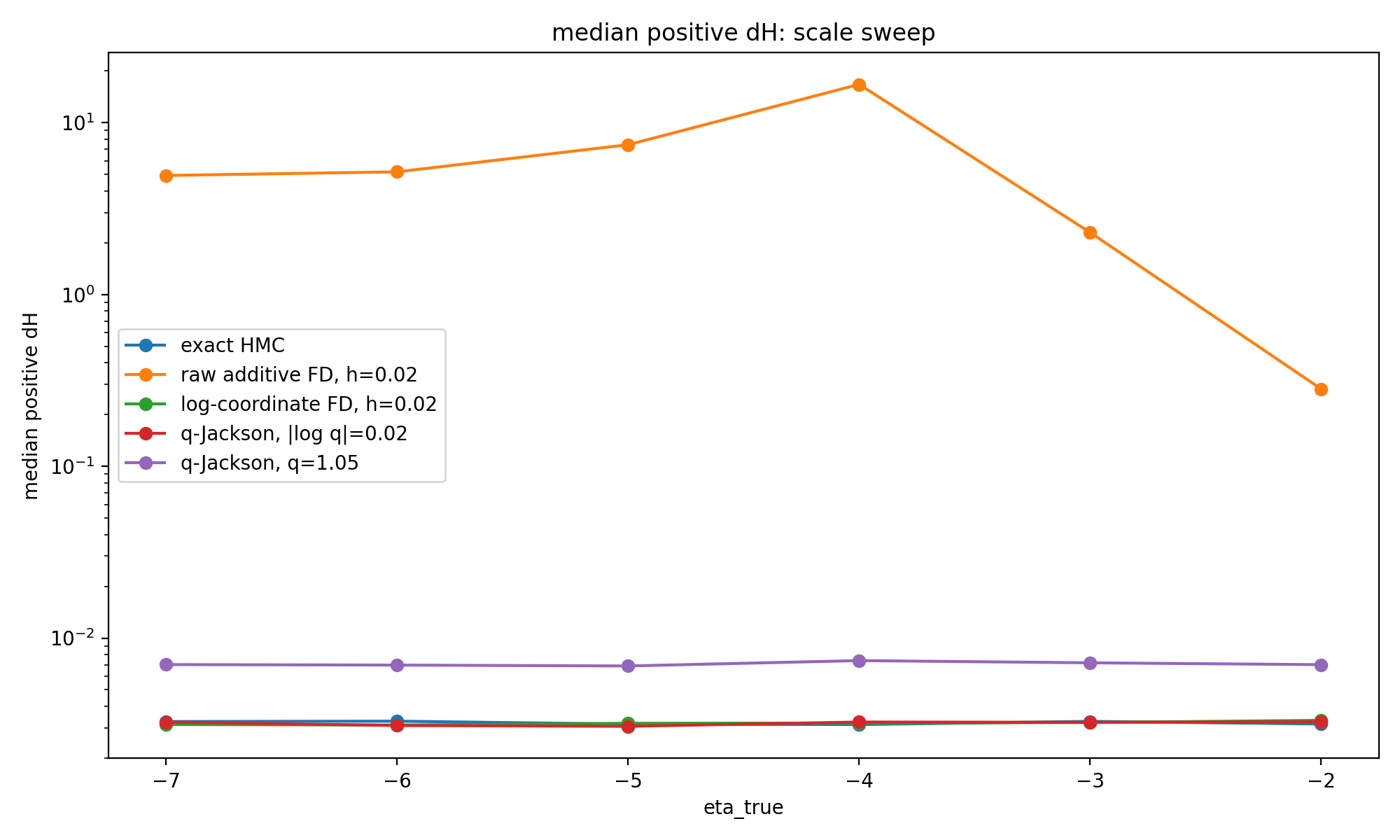}
        \caption{\(\mathrm{med.}(\Delta H_+)\)}
    \end{subfigure}

    \caption{
    Scale sweep for the variance-component positive-scale target with fixed
    dimension and varying true log scale \(\eta_{\rm true}\).
    }
    \label{fig:vc-scale-sweep}
\end{figure}

Figure~\ref{fig:vc-scale-sweep} reports the scale sweep.  The raw additive
finite-difference method is highly sensitive to the absolute scale of
\(\tau_i\).  When the posterior is concentrated at a small scale, a fixed
additive perturbation in \(\tau_i\) is no longer local relative to the posterior
scale.  This leads to large force errors, reduced acceptance rates, inflated
Hamiltonian errors, and degraded marginal accuracy.  By contrast, the
log-coordinate finite difference and the \(q\)-Jackson force remain stable
across the scale range because their perturbations are multiplicative rather
than additive.

\begin{figure}[p]
    \centering
    \begin{subfigure}{0.47\textwidth}
        \centering
        \includegraphics[width=\textwidth,height=0.23\textheight,keepaspectratio]{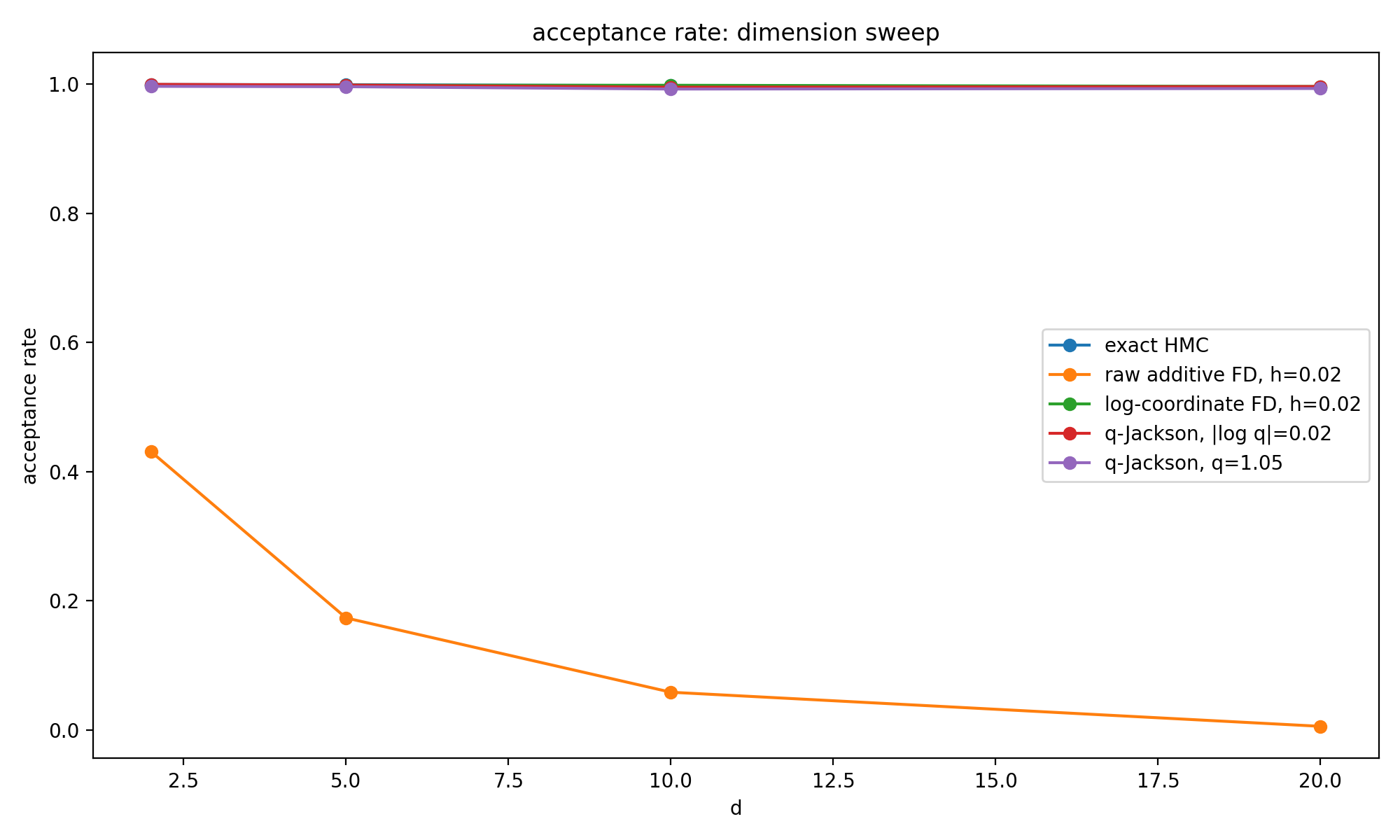}
        \caption{Acceptance rate}
    \end{subfigure}
    \hfill
    \begin{subfigure}{0.47\textwidth}
        \centering
        \includegraphics[width=\textwidth,height=0.23\textheight,keepaspectratio]{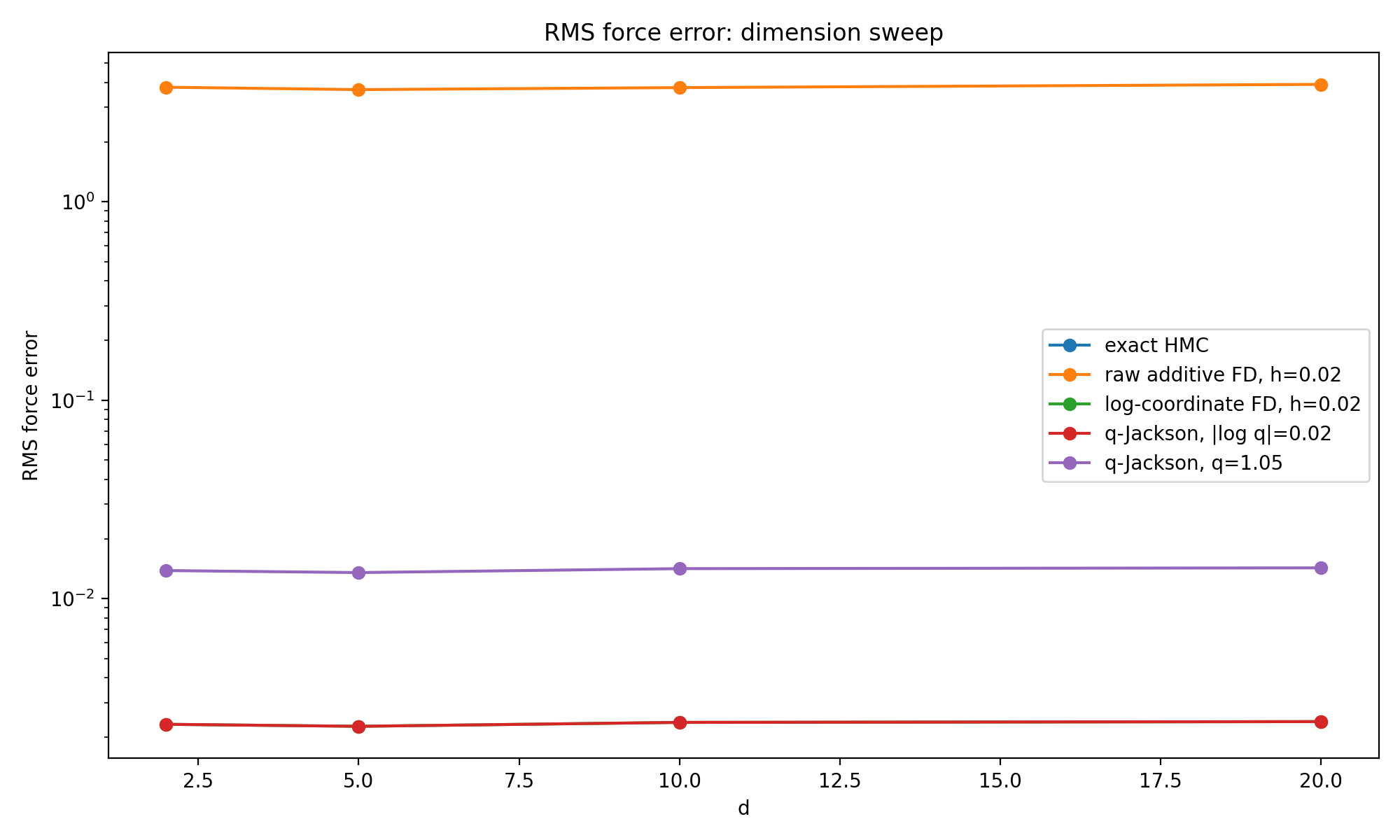}
        \caption{RMS force error}
    \end{subfigure}

    \vspace{0.25em}

    \begin{subfigure}{0.47\textwidth}
        \centering
        \includegraphics[width=\textwidth,height=0.23\textheight,keepaspectratio]{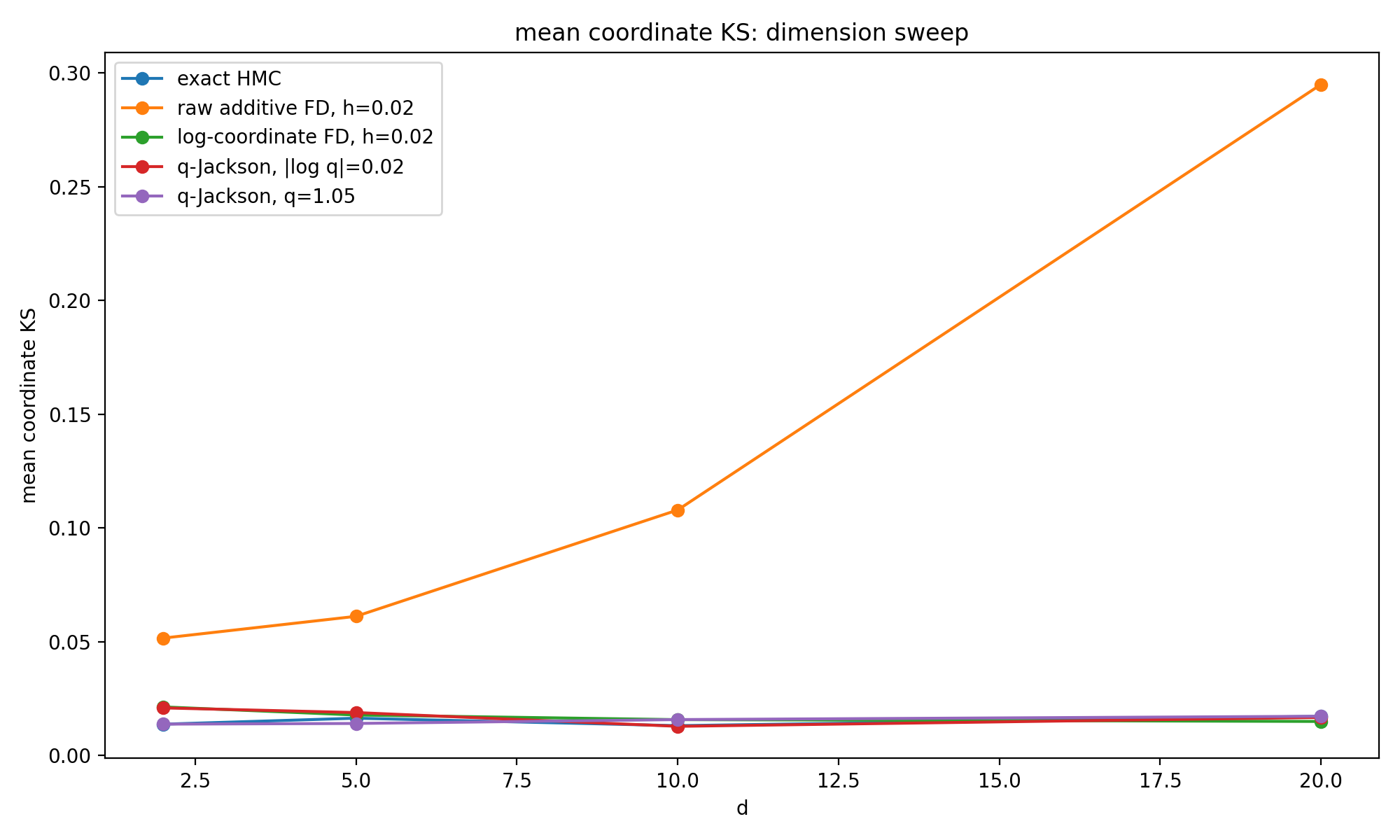}
        \caption{Mean coordinate KS distance}
    \end{subfigure}
    \hfill
    \begin{subfigure}{0.47\textwidth}
        \centering
        \includegraphics[width=\textwidth,height=0.23\textheight,keepaspectratio]{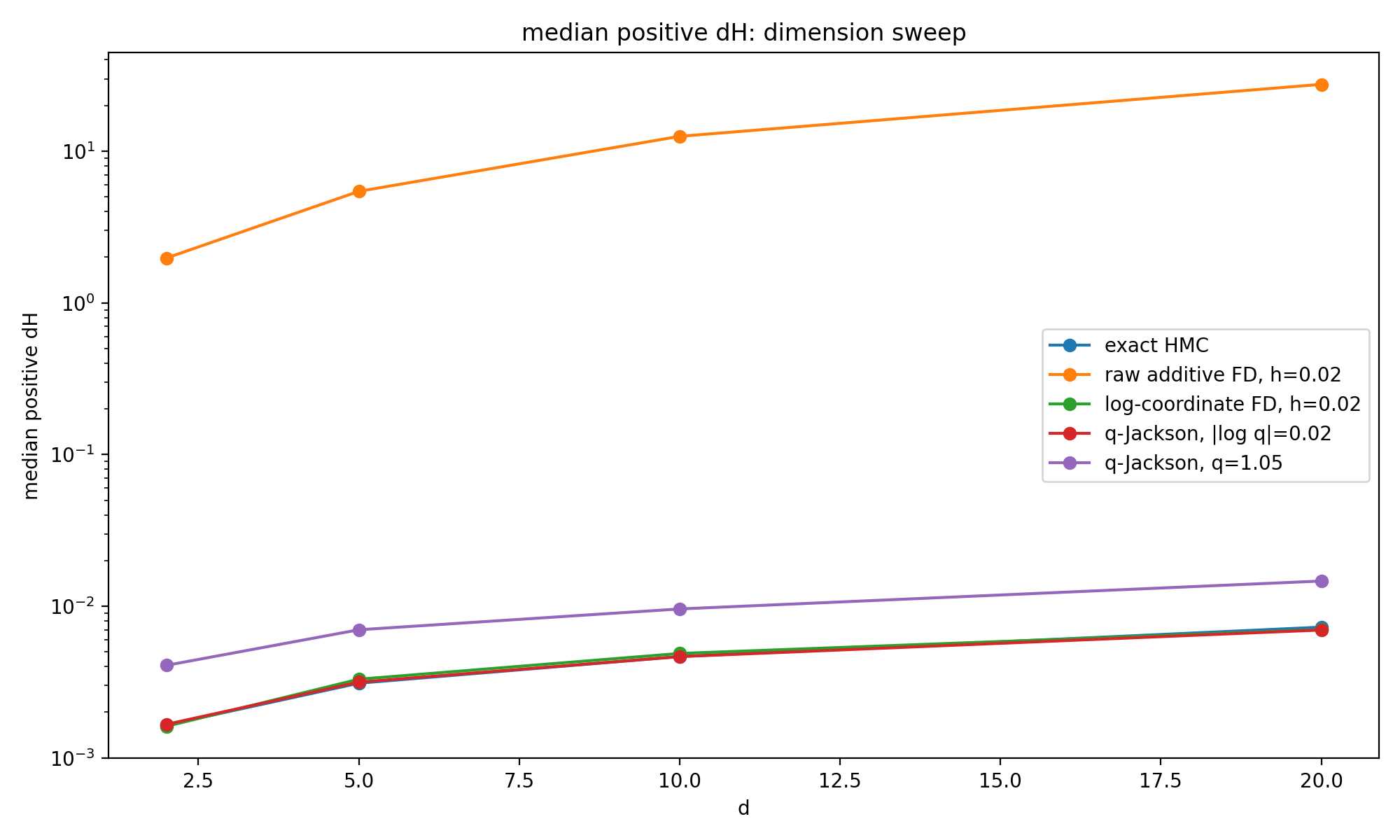}
        \caption{\(\mathrm{med.}(\Delta H_+)\)}
    \end{subfigure}

    \caption{
    Dimension sweep for the variance-component positive-scale target with fixed
    true log scale and increasing dimension \(d\).
    }
    \label{fig:vc-dim-sweep}
\end{figure}

Figure~\ref{fig:vc-dim-sweep} reports the dimension sweep.  The same qualitative
behavior persists as the number of variance components increases.  The raw
additive finite-difference force becomes increasingly unstable in both
acceptance and marginal accuracy, while the log-coordinate and \(q\)-Jackson
forces remain close to the exact-gradient reference.  This suggests that the
scale compatibility of the \(q\)-Jackson perturbation is not a one-dimensional
artifact.

Table~\ref{tab:hmc_comparison} reports two representative configurations: one
from the small-scale sweep and one from the high-dimensional sweep.  These
values provide a numerical summary of the trends in
Figures~\ref{fig:vc-scale-sweep}--\ref{fig:vc-dim-sweep}.  The key observation
is that the raw additive finite-difference method fails because its perturbation
size is not adapted to the scale of the positive parameter.  The \(q\)-Jackson
force avoids this problem by using relative perturbations.

\begin{table*}[htbp]
\centering
\small
\caption{
Comparison of HMC variants for the variance-component positive-scale
black-box target under scale sweep \((d=5,\eta_{\rm true}=-4)\) and dimension
sweep \((d=20,\eta_{\rm true}=-6)\). Smaller KS distance, force error, and
\(\mathrm{med.}(\Delta H_+)\) are better; larger acceptance rate and ESS are
better.
}
\label{tab:hmc_comparison}
\begin{tabular}{@{}lccccc@{}}
\toprule
\multirow{2}{*}{Method} & \multicolumn{5}{c}{Value} \\
\cmidrule(lr){2-6}
 & Acc. & KS & ESS & Force err. & \(\mathrm{med.}(\Delta H_+)\) \\
\midrule
\multicolumn{6}{@{}l}{\textit{Scale sweep: } \(d=5,\ \eta_{\rm true}=-4\)} \\
\midrule
exact HMC
& 0.998 & 0.015 & 2447 & 0
& \(3.14\times 10^{-3}\) \\
raw additive FD
& 0.027 & 0.354 & -- & \(1.01\times 10^{5}\)
& \(1.67\times 10^{1}\) \\
log-coordinate FD
& 0.999 & 0.016 & 2274 & \(2.28\times 10^{-3}\)
& \(3.20\times 10^{-3}\) \\
\(q\)-Jackson, \(|\log q|=0.02\)
& 0.998 & 0.017 & 2383 & \(2.28\times 10^{-3}\)
& \(3.25\times 10^{-3}\) \\
\(q\)-Jackson, \(q=1.05\)
& 0.995 & 0.017 & 2479 & \(1.36\times 10^{-2}\)
& \(7.39\times 10^{-3}\) \\
\midrule
\multicolumn{6}{@{}l}{\textit{Dimension sweep: } \(d=20,\ \eta_{\rm true}=-6\)} \\
\midrule
exact HMC
& 0.994 & 0.017 & 2484 & 0
& \(7.27\times 10^{-3}\) \\
raw additive FD
& 0.006 & 0.295 & -- & \(3.91\times 10^{0}\)
& \(2.76\times 10^{1}\) \\
log-coordinate FD
& 0.996 & 0.015 & 2488 & \(2.40\times 10^{-3}\)
& \(7.01\times 10^{-3}\) \\
\(q\)-Jackson, \(|\log q|=0.02\)
& 0.996 & 0.017 & 2488 & \(2.40\times 10^{-3}\)
& \(6.97\times 10^{-3}\) \\
\(q\)-Jackson, \(q=1.05\)
& 0.993 & 0.017 & 2428 & \(1.43\times 10^{-2}\)
& \(1.47\times 10^{-2}\) \\
\bottomrule
\end{tabular}
\end{table*}

Overall, this experiment identifies a regime in which the \(q\)-Jackson
construction has a clear computational interpretation.  For positive-scale
black-box targets, the multiplicative perturbation \(\tau_i\mapsto q\tau_i\)
naturally matches the geometry of scale parameters.  The resulting force
approximation behaves like a log-coordinate finite difference while remaining
expressible directly in the original positive variable.  The experiment
therefore supports the view that \(q\)-HMC is particularly useful when the
target is evaluated as a black-box function of positive scale parameters and
ordinary additive finite differences are poorly matched to the local scale of
the posterior.

\subsection{A mixed log-trap target}
\label{subsec:mixed-log-trap}

The preceding positive-scale examples show that the \(q\)-Jackson force has a
scale-consistent interpretation in logarithmic coordinates.  We now consider a
mixed target in which such a global logarithmic reparametrization is not
available.  The purpose of this example is to test whether the proposed
\(q\)-Jackson realization can still provide a stable proposal mechanism in a
coupled coordinate system containing both small positive variables and signed
variables.

Let
\[
x=(\theta,w)\in(0,\infty)^m\times\mathbb R^m ,
\]
where \(\theta_i>0\) are small positive variables and \(w_i\) are signed
variables.  We introduce normalized variables
\[
z_i=\frac{\theta_i}{s_i^\theta},\qquad
v_i=\frac{w_i}{s_i^w},
\]
where \(s_i^\theta>0\) and \(s_i^w>0\) are coordinate scales.  For
\(i=1,\ldots,m-2\), define the residual
\begin{align}
\label{eq:mixed-log-trap-residual}
r_i(\theta,w)
&=
\frac{v_i}{z_i+c_0}
+a_1 z_i
+a_2 v_{i+1}
+a_3 z_{i+1}v_i \notag\\
&\quad
+\eta\sin\left(
\frac{v_{i+1}}{z_i+\lambda z_{i+1}+c_0}
\right)
+a_4 z_{i+2}v_i v_{i+1}-d_i .
\end{align}
The target potential is
\begin{align}
\label{eq:mixed-log-trap-potential}
U(\theta,w)
=\frac{1}{2\sigma^2}
\sum_{i=1}^{m-2} r_i(\theta,w)^2+\frac{1}{2\tau_\theta^2}\sum_{i=1}^m (z_i-z_i^\dagger)^2+\frac{1}{2\tau_w^2}\sum_{i=1}^m (v_i-v_i^\dagger)^2 .
\end{align}
The data \(d_i\) are generated from \eqref{eq:mixed-log-trap-residual} at a
reference point and then perturbed by Gaussian noise.  In the experiment we set
\(m=12\), so that \(\dim(x)=24\).  The positive variables satisfy
\[
\theta_i\in[1.504\times10^{-8},\, 3.578\times10^{-4}],
\]
whereas the signed variables satisfy
\[
w_i\in[-6.343\times10^{-1},\, 1.620\times10^{-1}].
\]
Thus a global logarithmic transformation of the full state is not available.
A partial logarithmic transformation of the positive block is possible, but it
does not decouple the target and still leaves a multiscale signed block.

We compare exact-gradient HMC, fixed additive finite-difference HMC,
partial-log finite-difference HMC, scaled finite-difference HMC, and the
proposed \(q\)-Jackson HMC.  The partial-log method applies
\(\theta_i\mapsto \theta_i e^{\pm\varepsilon_\theta}\) to the positive block
and an additive perturbation \(w_j\mapsto w_j\pm h_w\) to the signed block.
The scaled finite difference uses a coordinatewise relative step.  The
\(q\)-Jackson force applies the symmetric multiplicative perturbation
\(x_i\mapsto e^{\pm|\log q|/2}x_i\) directly in the original mixed coordinates.
For \(\theta_i>0\), this perturbation preserves positivity; for nonzero signed
variables, it preserves the sign.

The scaled finite-difference baseline uses a coordinatewise relative additive
step
\[
    h_i=\varepsilon\max\{|x_i|,h_{\min}\},
\]
and approximates the potential gradient by
\[
    \widehat G_i^{\rm scaled}(x)
    =
    \frac{U(x+h_i e_i)-U(x-h_i e_i)}{2h_i}.
\]
This baseline is included to separate the effect of relative scaling from the
specific \(q\)-Jackson construction.

\begin{table*}[htbp]
\centering
\caption{
Results for the mixed log-trap target.  The force error is the median relative
error with respect to the exact gradient over test points.  Smaller force error
and \(\mathrm{med.}(\Delta H_+)\), and larger acceptance rate and ESS are preferred.
The column ``bad rej.'' records proposals rejected because the numerical
trajectory leaves the effective domain or gives a non-finite value.
}
\label{tab:mixed-log-trap}
\scriptsize
\setlength{\tabcolsep}{3.5pt}
\begin{tabular}{lccccc}
\toprule
Method & Acc. & ESS & Force err. & \(\mathrm{med.}(\Delta H_+)\) & bad rej. \\
\midrule
Exact-gradient HMC
& 0.694 & 29.48 & 0 & \(4.564{\times}10^{-1}\) & 8 \\
Additive FD, \(h=10^{-8}\)
& 0.000 & 1.00 & \(8.859{\times}10^{0}\) & \(5.086{\times}10^{1}\) & 0 \\
Additive FD, \(h=10^{-12}\)
& 0.684 & 29.70 & \(2.759{\times}10^{-8}\) & \(5.155{\times}10^{-1}\) & 1 \\
Partial-log FD, \(h_w=10^{-8}\)
& 0.636 & 19.91 & \(2.158{\times}10^{-3}\) & \(5.892{\times}10^{-1}\) & 57 \\
Partial-log FD, \(h_w=10^{-12}\)
& 0.435 & 17.07 & \(2.158{\times}10^{-3}\) & \(3.176{\times}10^{0}\) & 70 \\
Scaled FD, \(\varepsilon=0.02\)
& 0.742 & 35.70 & \(3.053{\times}10^{-3}\) & \(4.427{\times}10^{-1}\) & 8 \\
\(q\)-Jackson, \(|\log q|=0.02\)
& 0.695 & 35.06 & \(6.812{\times}10^{-4}\) & \(5.033{\times}10^{-1}\) & 7 \\
\(q\)-Jackson, \(|\log q|=0.005\)
& 0.705 & 36.20 & \(4.257{\times}10^{-5}\) & \(4.454{\times}10^{-1}\) & 25 \\
\bottomrule
\end{tabular}
\end{table*}

Table~\ref{tab:mixed-log-trap} shows three distinct effects.  First, a fixed
additive finite-difference step can be severely mismatched to the coordinate
scales.  With \(h=10^{-8}\), the median force error is \(8.859\) and the HMC
acceptance rate drops to zero.  A much smaller step \(h=10^{-12}\) gives a
highly accurate force approximation and restores the behavior of exact HMC, but
this step is problem-dependent and requires prior knowledge of the small
coordinate scales.

Second, a partial logarithmic treatment of the positive variables does not
remove the difficulty in the mixed coordinate system.  Although the force error
is reduced to \(2.158\times10^{-3}\), the method still produces many bad
rejections, and reducing the additive step in the signed block from \(10^{-8}\)
to \(10^{-12}\) does not improve the force error.  This indicates that the
remaining error is associated with the coupled mixed-coordinate structure, not
only with the finite-difference step used for the signed block.

Third, the \(q\)-Jackson force remains stable when applied directly in the
original coordinates.  For \(|\log q|=0.02\), the force error is
\(6.812\times10^{-4}\), and the resulting HMC transition has an acceptance rate
and ESS close to the exact-gradient benchmark.  Reducing the deformation level
to \(|\log q|=0.005\) further lowers the force error to
\(4.257\times10^{-5}\), while preserving comparable sampling behavior.  These
results show that the \(q\)-Jackson realization provides a relative
perturbation mechanism that is not tied to a global logarithmic
reparametrization.

This experiment should not be interpreted as showing that the \(q\)-Jackson
force uniformly dominates all finite-difference approximations.  A carefully
tuned additive finite-difference step can be very accurate.  Rather, the point
is that the \(q\)-deformed construction naturally supplies a relative
perturbation rule in the original mixed coordinates.  This is useful when fixed
additive perturbations suffer from scale mismatch, a global logarithmic
transformation is unavailable, and a partial logarithmic transformation does
not decouple the target.


\subsection{Positive-coefficient diffusion inverse problem}
\label{subsec:inverse-positive-coefficient}

We now turn to a PDE-constrained Bayesian inverse problem with a strictly positive functional parameter. This example serves two purposes.
 First, it clarifies
how the multiplicative perturbation underlying the \(q\)-Jackson force is
connected with the standard adjoint sensitivity calculus for positive
coefficients.  Second, it provides a computational-physics benchmark in which
the potential is generated by repeated PDE solves rather than by an explicit
closed-form density.

We begin with a general formulation.  Let \(u\) be the state variable and let
\(\alpha\) be a positive system parameter appearing in the PDE constraint
\begin{align}
    \mathcal F(u;\alpha)=0 .
\end{align}
Given observations \(\mathbf d_j\) at locations \(x_j\), consider the data
misfit
\begin{align}
    \mathcal J(\alpha)
    =
    \frac12\sum_j |u(x_j;\alpha)-\mathbf d_j|^2 .
\end{align}
For a fixed point \(\xi\in\Omega\), the formal \(q\)-Jackson perturbation of
the coefficient is written as
\begin{align}
    \alpha_{q,\xi}(x)
    =
    \begin{cases}
        \alpha(x), & x\neq \xi,\\
        q^2\alpha(\xi), & x=\xi .
    \end{cases}
\end{align}
The corresponding \(q\)-Jackson sensitivity is
\begin{align}
\label{eq:q-jackson-functional-sensitivity}
    \mathcal D_{\alpha(\xi)}\mathcal J(\alpha)
    :=
    \frac{\mathcal J(\alpha_{q,\xi})-\mathcal J(\alpha)}
    {(q^2-1)\alpha(\xi)} .
\end{align}
The pointwise notation should be understood either after spatial discretization
or as the limit of a localized perturbation.  In a finite-dimensional
discretization, \eqref{eq:q-jackson-functional-sensitivity} simply corresponds
to multiplying one positive coefficient degree of freedom by \(q^2\).

Let \(\delta u_{q,\xi}=u(\alpha_{q,\xi})-u(\alpha)\).  Linearizing the forward
equation at \((u,\alpha)\) gives
\begin{align}
    L\delta u_{q,\xi}
    =
    -(q^2-1)\alpha(\xi)A_\xi ,
\end{align}
where
\[
    L=\frac{\partial\mathcal F}{\partial u},
    \qquad
    A_\xi=\frac{\partial\mathcal F}{\partial\alpha}\delta_\xi .
\]
Let \(\lambda\) solve the adjoint equation
\begin{align}
    L^*\lambda
    =
    \sum_j
    \bigl(u(x_j)-\mathbf d_j\bigr)\delta_{x_j}.
\end{align}
Then, up to higher-order terms,
\begin{align}
\begin{aligned}
    &\mathcal J(\alpha_{q,\xi})-\mathcal J(\alpha)\\
    &=\sum_j\bigl(u(x_j)-\mathbf d_j\bigr)
    \delta u_{q,\xi}(x_j)+O(\|\delta u_{q,\xi}\|^2)  \\
    &=
    \langle L^*\lambda,\delta u_{q,\xi}\rangle+O(\|\delta u_{q,\xi}\|^2) \\
    &=
    \langle \lambda,L\delta u_{q,\xi}\rangle+O(\|\delta u_{q,\xi}\|^2) \\
    &=
    -(q^2-1)\alpha(\xi)
    \left[\left(\frac{\partial\mathcal F}{\partial\alpha}\right)^*
    \lambda\right](\xi)+O(\|\delta u_{q,\xi}\|^2).
    \end{aligned}
\end{align}
Therefore, to first order,
\begin{align}
\label{eq:q-jackson-adjoint-consistency}
    \mathcal D_{\alpha(\xi)}\mathcal J(\alpha)
    =
    -\mathcal G_{\alpha}^*(\lambda)(\xi),
\end{align}
where
\[
    \mathcal G_{\alpha}^*(\lambda)
    :=
    \left(\frac{\partial\mathcal F}{\partial\alpha}\right)^*
    \lambda .
\]
Thus the \(q\)-Jackson sensitivity coincides, at the first-order level, with
the classical adjoint sensitivity in the small-perturbation limit.  This
calculation is not meant to replace the adjoint method.  Rather, it explains
why a multiplicative \(q\)-Jackson perturbation is a natural black-box
perturbation rule for positive functional parameters.

We now specialize this observation to the diffusion coefficient inverse
problem
\begin{align}
\label{eq:diffusion-forward}
\begin{cases}
    -\dfrac{d}{dx}\left(\alpha(x)\dfrac{du}{dx}\right)=f(x),
    & x\in(0,1),\\[4pt]
    u(0)=u(1)=0,
\end{cases}
\end{align}
where \(\alpha(x)>0\).  Noisy observations are taken at sensor locations
\(x_j\):
\[
    d_j=u(x_j;\alpha)+\eta_j,
    \qquad
    \eta_j\sim N(0,\sigma_{\rm obs}^2).
\]
The data misfit is
\begin{align}
\label{eq:diffusion-misfit}
    \Phi(\alpha)
    =
    \frac{1}{2\sigma_{\rm obs}^2}
    \sum_j
    \left|u(x_j;\alpha)-d_j\right|^2 .
\end{align}
For the continuous diffusion problem, the adjoint variable satisfies
\begin{align}
\begin{cases}
    -\nabla\cdot(\alpha\nabla\lambda)
    =
    \displaystyle
    \frac{1}{\sigma_{\rm obs}^2}
    \sum_j
    \bigl(u(x_j)-d_j\bigr)\delta_{x_j},
    & x\in\Omega,\\
    \lambda=0, & x\in\partial\Omega .
\end{cases}
\end{align}
The corresponding coefficient sensitivity is
\begin{align}
    \frac{\delta \Phi}{\delta \alpha}(\xi)
    =
    -\nabla u(\xi)\cdot\nabla\lambda(\xi).
\end{align}
Hence, in this representative positive-coefficient inverse problem, the
formal \(q\)-Jackson sensitivity recovers the usual adjoint sensitivity in the
small-perturbation limit.

In the numerical implementation, the coefficient is represented by positive
finite-volume edge coefficients
\[
    \alpha=(\alpha_0,\ldots,\alpha_n),
    \qquad
    \alpha_e>0,
\]
and the sampler is formulated in the log-coordinate
\[
    y_e=\log\alpha_e,
    \qquad e=0,\ldots,n .
\]
Let \(h=1/(n+1)\), and let \(u_i\) denote the discrete state at the interior
nodes \(x_i=ih\), \(i=1,\ldots,n\).  With \(u_0=u_{n+1}=0\), the finite-volume
discretization of \eqref{eq:diffusion-forward} gives
\begin{align}
\label{eq:discrete-forward}
\begin{aligned}
   & (K(\alpha)u)_i = \frac{
    \alpha_{i-1}(u_i-u_{i-1}) + \alpha_i(u_i-u_{i+1})}{h^2}=f_i, 
    \\ 
    &i=1,\ldots, n .
    \end{aligned}
\end{align}
Let \(C\) be the observation matrix that extracts the observed components of
\(u\).  The discrete data misfit is
\begin{align}
\label{eq:discrete-misfit}
    \Phi_h(\alpha)
    =
    \frac{1}{2\sigma_{\rm obs}^2}
    \|Cu(\alpha)-d\|^2 .
\end{align}
We impose a Gaussian prior on the log-coefficient,
\[
    y\sim N(m,\sigma_a^2 I).
\]
The posterior potential with respect to the log-coordinate \(y\) is therefore
\begin{align}
\label{eq:diffusion-log-potential}
    V(y)
    =
    \Phi_h(\exp y)
    +
    \frac{1}{2\sigma_a^2}\|y-m\|^2 .
\end{align}
The HMC sampler is applied to \(V(y)\).  To express black-box perturbations
directly in the original positive coefficient variable \(\alpha\), we also use
the potential with respect to the measure \(d\alpha\):
\begin{align}
\label{eq:diffusion-alpha-potential}
    U_\alpha(\alpha)
    =
    V(\log\alpha)
    +
    \sum_{e=0}^n \log\alpha_e .
\end{align}
The additional term in \eqref{eq:diffusion-alpha-potential} is the Jacobian
contribution of the transformation \(\alpha=\exp y\).

The discrete adjoint derivative provides an exact-gradient benchmark.  Given
\(\alpha\), the forward state is obtained from
\[
    K(\alpha)u=f .
\]
The adjoint variable \(\lambda\) solves
\begin{align}
\label{eq:discrete-adjoint}
    K(\alpha)^T\lambda
    =
    \frac{1}{\sigma_{\rm obs}^2}
    C^T(Cu-d).
\end{align}
Then the derivative of the data misfit with respect to the edge coefficient
\(\alpha_e\) is
\begin{align}
\label{eq:adjoint-alpha-gradient}
    \frac{\partial \Phi_h}{\partial \alpha_e}
    =
    -\lambda^T
    \frac{\partial K(\alpha)}{\partial \alpha_e}
    u .
\end{align}
Equivalently, setting \(u_0=u_{n+1}=0\) and
\(\lambda_0=\lambda_{n+1}=0\), one obtains
\begin{align}
\label{eq:adjoint-edge-gradient}
    \frac{\partial \Phi_h}{\partial \alpha_e}
    =
    -
    \frac{(\lambda_e-\lambda_{e+1})(u_e-u_{e+1})}{h^2},
    \qquad e=0,\ldots,n .
\end{align}
By the chain rule, the exact force in the log-coordinate is
\begin{align}
\label{eq:adjoint-log-force}
    F^{\rm adj}_e(y)=-\frac{\partial V}{\partial y_e}
    = -\left[
    \alpha_e
    \frac{\partial \Phi_h}{\partial \alpha_e}+\frac{y_e-m_e}{\sigma_a^2}
    \right],
    \,
    \alpha_e=e^{y_e}.
\end{align}
This is the reference force used in the exact-gradient HMC benchmark.

For a positive coefficient degree of freedom \(\alpha_e\), the symmetric
\(q\)-Jackson perturbation used in the numerical sampler is multiplicative:
\[
    \alpha_e\mapsto q\alpha_e,
    \qquad
    \alpha_e\mapsto \alpha_e/q .
\]
Let
\[
    \alpha^{(e,+)}
    =
    (\alpha_0,\ldots,q\alpha_e,\ldots,\alpha_n),
    \,\,
    \alpha^{(e,-)}
    =
    (\alpha_0,\ldots,\alpha_e/q,\ldots,\alpha_n).
\]
The corresponding symmetric \(q\)-Jackson sensitivity of the data misfit is
\begin{align}
\label{eq:q-jackson-misfit-sensitivity}
    D^{\rm sym}_{q,\alpha_e}\Phi_h(\alpha)
    =
    \frac{
    \Phi_h(\alpha^{(e,+)})
    -
    \Phi_h(\alpha^{(e,-)})
    }{
    2\alpha_e\log q
    } .
\end{align}
Since this is a centered finite difference in the logarithmic coordinate
\(\log\alpha_e\), Taylor expansion gives
\begin{align}
\label{eq:q-jackson-adjoint-consistency-discrete}
    D^{\rm sym}_{q,\alpha_e}\Phi_h(\alpha)
    =
    \frac{\partial \Phi_h}{\partial \alpha_e}
    +
    O\!\left((\log q)^2\right),
    \qquad q\to1 .
\end{align}
Combining \eqref{eq:adjoint-alpha-gradient} and
\eqref{eq:q-jackson-adjoint-consistency-discrete}, the \(q\)-Jackson
sensitivity is consistent with the discrete adjoint sensitivity in the
small-perturbation limit:
\begin{align}
\label{eq:q-jackson-adjoint-consistency-final}
    D^{\rm sym}_{q,\alpha_e}\Phi_h(\alpha)
    =
    -\lambda^T
    \frac{\partial K(\alpha)}{\partial \alpha_e}
    u
    +
    O\!\left((\log q)^2\right).
\end{align}

For sampling, the force must approximate \(-\nabla_y V(y)\).  When the
black-box potential is evaluated in the original coefficient variable through
\(U_\alpha(\alpha)\), the \(q\)-Jackson approximation to
\(\partial V/\partial y_e\) is
\begin{align}
\label{eq:q-jackson-posterior-force}
    \widehat{\nabla_y V}^{\,q}_e(y)
    =
    \frac{
    U_\alpha(\alpha^{(e,+)})
    -
    U_\alpha(\alpha^{(e,-)})
    }{2\log q}
    -1 .
\end{align}
The subtraction of \(1\) removes the derivative of the Jacobian term
\(\sum_e\log\alpha_e\) in \eqref{eq:diffusion-alpha-potential}.  The
corresponding \(q\)-Jackson force is
\[
    F^q_e(y)
    =
    -\widehat{\nabla_y V}^{\,q}_e(y).
\]
As \(q\to1\), this force converges to the exact log-coordinate force
\eqref{eq:adjoint-log-force} with second-order error in \(\log q\), provided
the discretized forward map is sufficiently smooth.

We compare five force constructions.  The exact-gradient benchmark uses the
adjoint force \eqref{eq:adjoint-log-force}, which requires one forward solve
and one adjoint solve per force evaluation.  The raw additive finite-difference
force perturbs the original coefficient additively,
\[
    \alpha_e\mapsto \alpha_e\pm h_{\rm abs}.
\]
When the backward perturbation violates positivity, a one-sided difference is
used.  This method is included as a diagnostic for scale mismatch in the
original positive variable.  The log-coordinate finite-difference force uses
\[
    y_e\mapsto y_e\pm h_{\log}.
\]
Finally, the \(q\)-Jackson force uses the multiplicative perturbation
\[
    \alpha_e\mapsto q\alpha_e,\,
    \alpha_e\mapsto \alpha_e/q .
\]
We report both a local perturbation with \(|\log q|=h_{\log}\) and a fixed
deformation value \(q=1.05\).

Synthetic data are generated from the smooth positive coefficient
\begin{align}
\label{eq:true-diffusion-coefficient}
\begin{aligned}
   & y^\dagger(x) =-5+0.55\sin(2\pi x)+0.25\cos(4\pi x),\\
    &\alpha^\dagger(x)=\exp(y^\dagger(x)).
    \end{aligned}
\end{align}
In the reported experiment, \(n=20\), so that the positive coefficient is
represented by \(21\) edge coefficients.  We use \(10\) observation locations,
\(\sigma_{\rm obs}=0.015\), and \(\sigma_a=0.7\).  The source term is set to
\[
    f_i=10^{-2},\qquad i=1,\ldots,n,
\]
which gives a more informative data regime than the smaller forcing case.  The
finite-difference and \(q\)-Jackson parameters are
\(h_{\rm abs}=h_{\log}=0.02\), with an additional fixed deformation case
\(q=1.05\).  Each method is run with three chains; the first \(200\) iterations
of each chain are discarded as burn-in.

For each method, we report the acceptance rate, the mean coordinate effective
sample size, the RMS force error relative to the adjoint force,
\begin{align}
\label{eq:diffusion-force-error}
    E_{\rm force}
    =
    \left[
    \frac{1}{M}
    \sum_{r=1}^M
    \frac{
    \|F^{\rm method}(y^{(r)})-F^{\rm adj}(y^{(r)})\|^2
    }{n+1}
    \right]^{1/2},
\end{align}
and the median positive Hamiltonian error,
\begin{align}
\label{eq:median-positive-dH}
    \mathrm{med.}(\Delta H_+)
    =
    {\rm median}\{\Delta H_r:\Delta H_r>0\}.
\end{align}
We also report reconstruction errors for the posterior mean, including the
relative coefficient error
\[
    \frac{\|\bar\alpha-\alpha^\dagger\|}{\|\alpha^\dagger\|}
\]
and the root mean square error of the log-coefficient.

\begin{table*}[htbp]
\centering
\caption{
Numerical results for the positive-coefficient diffusion inverse problem with
\(f_i=10^{-2}\). The force error is measured relative to the adjoint force.
Smaller force error, \(\mathrm{med.}(\Delta H_+)\), relative coefficient error, and
log-coefficient RMSE are preferred, while larger acceptance rate and ESS are
preferred.
}
\label{tab:diffusion-inverse-results}
\scriptsize
\setlength{\tabcolsep}{4pt}
\renewcommand{\arraystretch}{1.08}
\begin{tabular}{@{}lcccccc@{}}
\toprule
Method
& Acc.
& ESS
& Force
& \(\mathrm{med.}(\Delta H_+)\)
& Rel. err.
& \(y\)-RMSE \\
\midrule
Exact adjoint HMC
& 0.998
& 61.63
& 0
& \(1.09{\times}10^{-3}\)
& 0.291
& 0.307 \\
Raw additive FD, \(h=0.02\)
& 0.159
& 19.03
& \(2.34{\times}10^{1}\)
& \(8.78{\times}10^{0}\)
& 0.342
& 0.386 \\
Log-coordinate FD, \(h=0.02\)
& 0.999
& 60.83
& \(3.18{\times}10^{-3}\)
& \(1.21{\times}10^{-3}\)
& 0.259
& 0.277 \\
\(q\)-Jackson, \(|\log q|=0.02\)
& 1.000
& 62.97
& \(3.18{\times}10^{-3}\)
& \(1.07{\times}10^{-3}\)
& 0.312
& 0.306 \\
\(q\)-Jackson, \(q=1.05\)
& 0.997
& 62.13
& \(1.89{\times}10^{-2}\)
& \(4.39{\times}10^{-3}\)
& 0.287
& 0.286 \\
\bottomrule
\end{tabular}
\end{table*}

Table~\ref{tab:diffusion-inverse-results} reports the numerical results.  The
exact adjoint HMC provides the reference force.  The raw additive
finite-difference force has a substantially larger force error than all
multiplicative or log-coordinate approximations.  Its RMS force error is
\(2.34\times10^{1}\), whereas the log-coordinate finite difference and the
symmetric \(q\)-Jackson force with \(|\log q|=0.02\) both have force error
\(3.18\times10^{-3}\).  This confirms that, for a positive coefficient
parameter, the multiplicative \(q\)-Jackson perturbation reproduces the
log-coordinate sensitivity much more accurately than a raw additive
perturbation in the original coefficient variable.

The same effect is reflected in the Hamiltonian error.  The median positive
Hamiltonian error of the raw additive finite-difference method is
\(8.78\times10^{0}\), while the corresponding value for the \(q\)-Jackson force
with \(|\log q|=0.02\) is \(1.07\times10^{-3}\), essentially at the same scale
as the exact adjoint and log-coordinate finite-difference benchmarks.  The
large energy error of the raw additive method also leads to a much lower
acceptance rate, \(0.159\), whereas the adjoint, log-coordinate, and
\(q\)-Jackson variants all maintain acceptance rates close to one under the
same integration parameters.

The reconstruction errors are of comparable magnitude across the stable
methods.  In this test, the log-coordinate finite difference gives the smallest
relative coefficient error, while the two \(q\)-Jackson variants remain close
to the adjoint benchmark.  Therefore, the main conclusion is not that
\(q\)-Jackson improves the coefficient reconstruction relative to adjoint HMC.
Rather, the experiment shows that the \(q\)-Jackson force is a
scale-consistent black-box approximation to the adjoint force for positive
functional parameters, while raw additive finite differences can introduce a
severe scale mismatch.

\begin{figure}[htbp]
    \centering
    \includegraphics[width=0.48\textwidth]{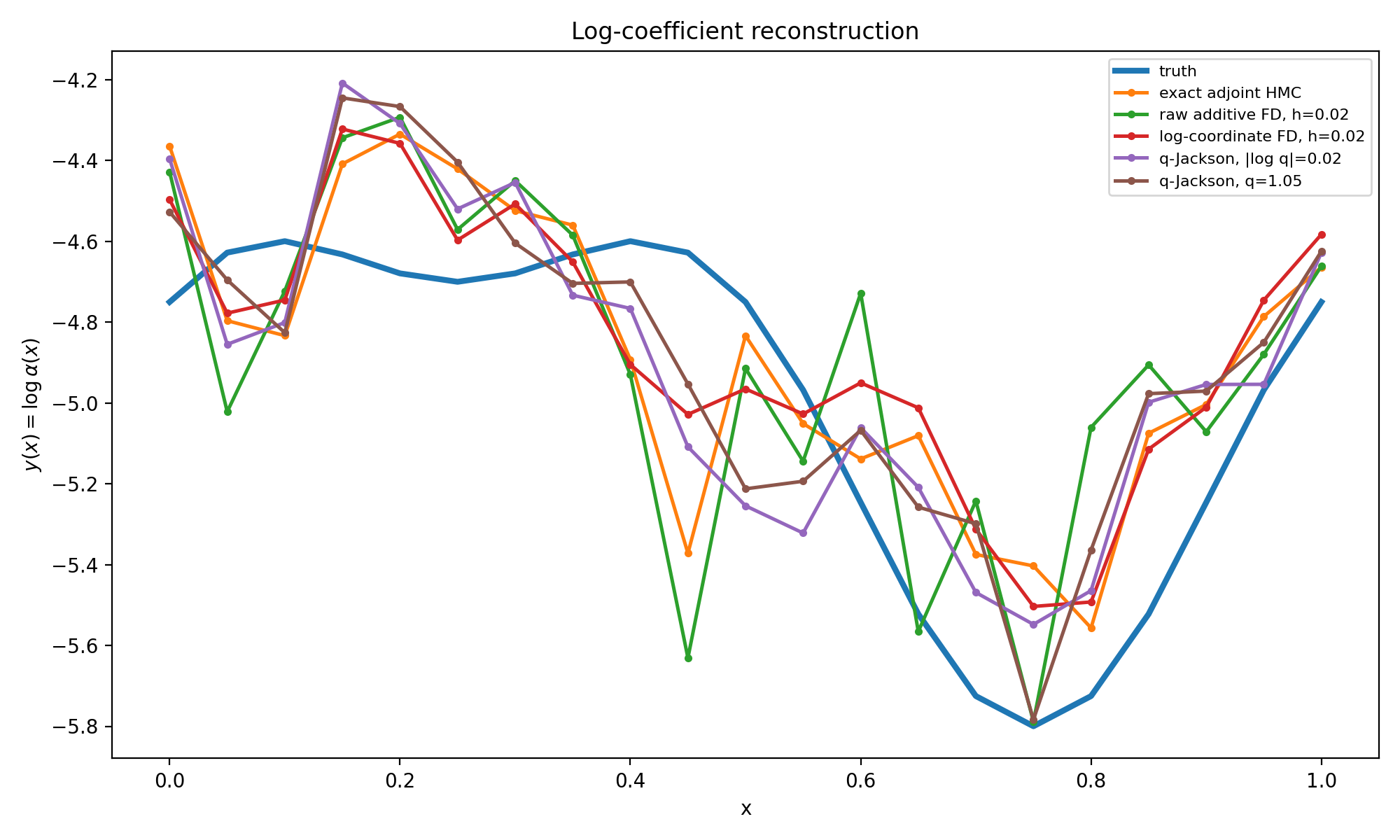}
    \includegraphics[width=0.48\textwidth]{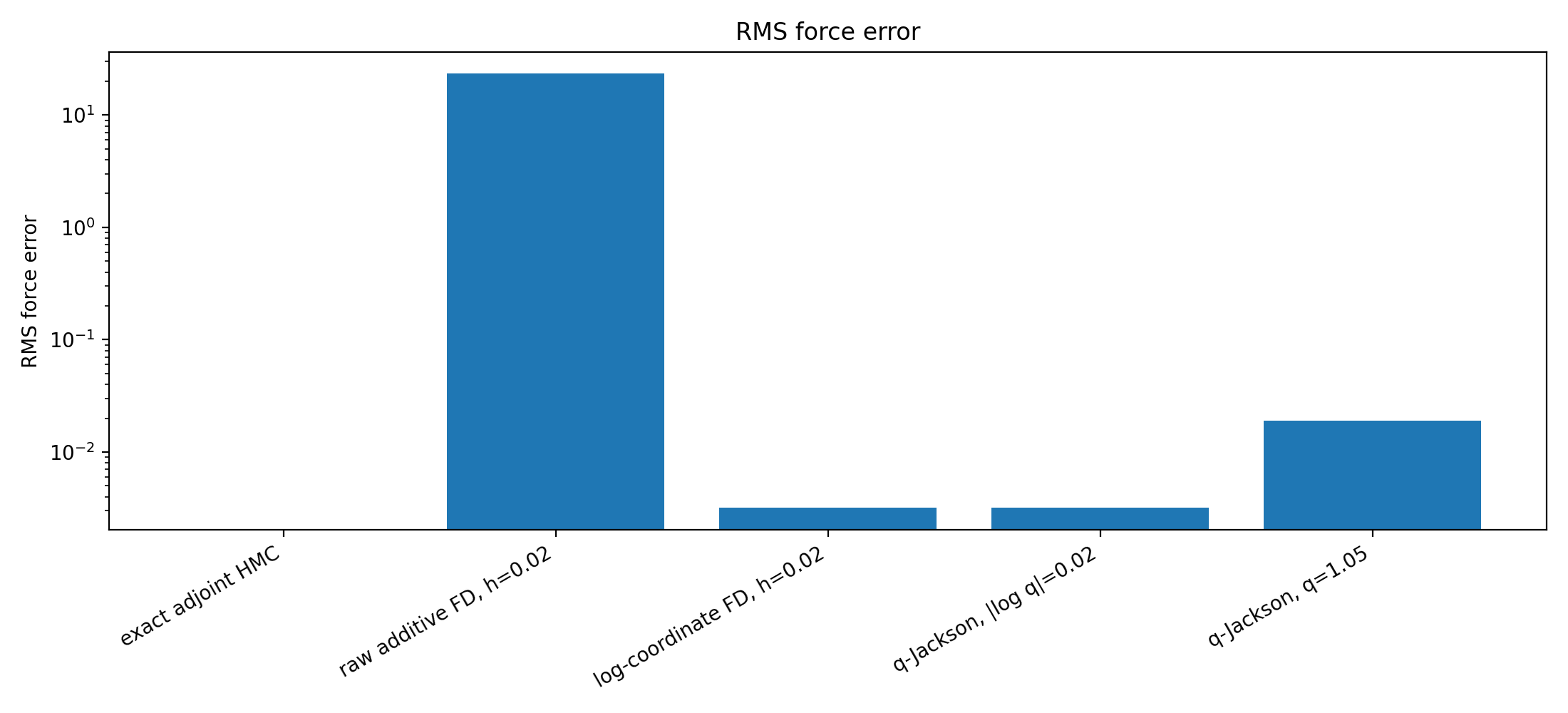}
    \caption{
    Positive-coefficient diffusion inverse problem.  Left: posterior mean
    reconstruction of the log-diffusion coefficient.  Right: RMS force error
    relative to the adjoint force.  The \(q\)-Jackson force uses multiplicative
    perturbations of the positive coefficient and closely matches the
    log-coordinate finite-difference force.
    }
    \label{fig:diffusion-reconstruction-force}
\end{figure}

\begin{figure}[htbp]
    \centering
    \includegraphics[width=0.55\textwidth]{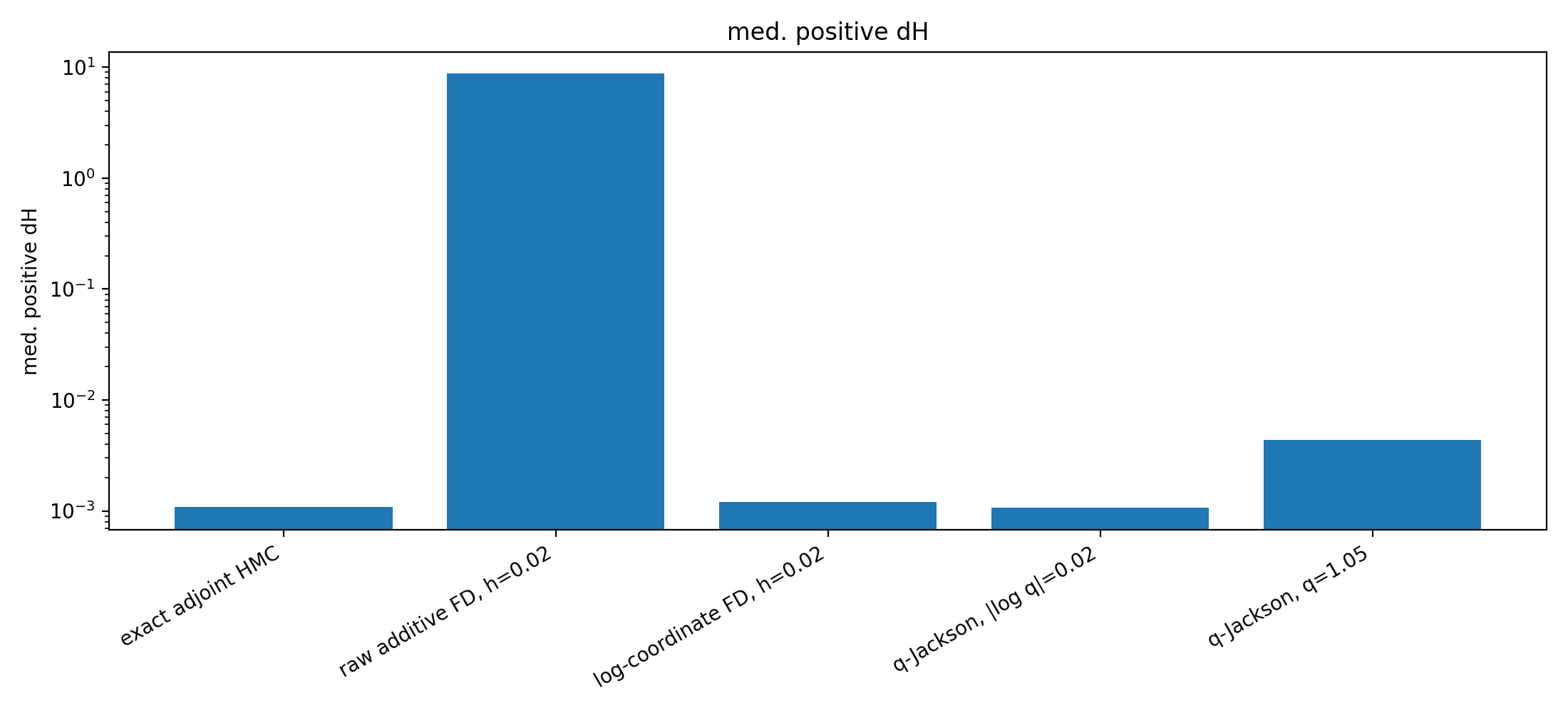}
    \caption{
    Median positive Hamiltonian error for the positive-coefficient diffusion
    inverse problem.  The raw additive finite-difference force produces a much
    larger energy error, whereas the \(q\)-Jackson force remains close to the
    adjoint and log-coordinate benchmarks.
    }
    \label{fig:diffusion-hamiltonian-error}
\end{figure}

This example should not be interpreted as replacing the adjoint method when an
adjoint solver is available.  Instead, it shows that the \(q\)-Jackson
multiplicative perturbation is consistent with adjoint sensitivity for positive
functional parameters and can therefore be used as a scale-compatible
black-box force approximation when only repeated evaluations of the
coefficient-to-misfit map are available.

\section{Conclusions}
\label{sec:conclude}

This paper developed a \(q\)-deformed Hamiltonian Monte Carlo framework
motivated by \(q\)-deformed phase-space mechanics. Starting from a Lagrangian
formulation, we derived the \(q\)-Hamiltonian equations and proved the formal
invariance of the associated \(q\)-symplectic form within the
\(q\)-deformed differential calculus. These results provide the geometric basis
for the proposed dynamics.
For practical sampling, the formal \(q\)-differential operators were realized
through Jackson-type differences. This leads to a computable
Metropolis-adjusted \(q\)-Hamiltonian proposal. The Metropolis correction, not
exact preservation of the formal \(q\)-geometric structure after numerical
realization, guarantees the validity of the Markov transition. We established
detailed balance and recorded sufficient nondegeneracy conditions for
irreducibility of the marginal position chain.

The numerical experiments clarify the computational role of the method. The
\(q\)-Jackson force provides a relative perturbation mechanism in the original
coordinates. For positive variables, it corresponds to a centered finite
difference in logarithmic coordinates; for mixed positive/signed targets, it can
be applied directly without requiring a global logarithmic reparametrization.
Across the tested examples, raw additive finite differences can suffer from
domain violations, scale mismatch, large force errors, and large Hamiltonian
errors. The \(q\)-Jackson force avoids these failures and tracks exact-gradient,
adjoint-force, log-coordinate, or scaled finite-difference benchmarks when such
references are available.

The proposed method is not intended as a universal replacement for
exact-gradient HMC, adjoint-based HMC, or modern adaptive HMC variants. Its main
contribution is to connect \(q\)-deformed Hamiltonian mechanics with a valid
Metropolis-corrected sampling scheme and to show that the induced Jackson-type
force can be useful for black-box or scale-sensitive targets. Future work should
study adaptive selection of the deformation parameter, integration with standard
HMC diagnostics, and applications to larger Bayesian models with constrained,
multiscale, or mixed-coordinate parameters.

\bibliography{sn-bibliography}


\end{document}